\documentclass{amsart}

\usepackage{amsthm,amsfonts,amsmath,amssymb,amscd,mathrsfs,verbatim,hyperref}

\numberwithin{equation}{section}
\newtheorem{theorem}{Theorem}[section]
\newtheorem{corollary}[theorem]{Corollary}
\newtheorem{proposition}[theorem]{Proposition}
\newtheorem{definition}{Definition}[section]

\newtheorem{lemma}[theorem]{Lemma}

\renewcommand{\And}{\qquad\text{and}\qquad}
\renewcommand{\and}{\quad\text{and}\quad}
\DeclareMathOperator{\MS}{\mathit{MS}}
\newcommand{\X}{\mathbb X}
\newcommand{\Y}{\mathbb Y}
\newcommand{\V}{\mathbb V}
\newcommand{\W}{\mathbb W}
\newcommand{\Z}{\mathbb Z}
\newcommand{\N}{\mathbb N}
\newcommand{\Q}{\mathbb Q}
\newcommand{\Complex}{\mathbb C}
\newcommand{\Symm}{\mathfrak{S}}
\newcommand{\abs}[1]{\lvert#1\rvert}
\newcommand{\la}{\lambda}
\newcommand{\qbin}[2]{\genfrac{[}{]}{0pt}{}{#1}{#2}}
\newcommand{\g}{\mathsf{g}}
\newcommand{\f}{\mathsf{f}}
\newcommand{\PP}{\mathsf{P}}
\newcommand{\QQ}{\mathsf{Q}}
\newcommand{\EE}{\mathsf{E}}
\newcommand{\FF}{\mathsf{F}}
\newcommand{\Skew}[4]{\EE_{#1/#2}(#3,#4)}
\newcommand{\SkewF}[4]{\FF_{#1/#2}(#3,#4)}
\newcommand{\Skewqt}[6]{\EE_{#1/#2}(#3,#4;#5,#6)}
\newcommand{\RR}{\mathsf{R}}
\newcommand{\MM}{\mathsf{M}}
\newcommand{\MMS}{\mathsf{MS}}
\renewcommand{\AA}{\mathbb{A}}
\newcommand{\BB}{\mathbb{B}}
\newcommand{\spec}[1]{\langle #1\rangle}
\newcommand{\gn}{\mathfrak{gl}_n}
\newcommand{\sln}{\mathfrak{sl}_n}
\newcommand*{\urlw}[1]{\href{http://www.#1}{\nolinkurl{www.#1}}}

\begin{document}

\title[Interpolation Macdonald polynomials]
{Nonsymmetric interpolation Macdonald polynomials and $\gn$ basic 
hypergeometric series}

\author{Alain Lascoux, Eric M. Rains and S. Ole Warnaar}
\thanks{Work supported by the ANR project MARS (BLAN06-2 134516), the
National Science Foundation grant DMS-0401387, and the 
Australian Research Council}

\address{CNRS, Institut Gaspard Monge, Universit\'{e} Paris-Est,
Marne-La-Vall\'ee, France}
\address{Department of Mathematics, California Institute of Technology,
Pasadena, CA 91125, USA} 
\address{Department of Mathematics, The University of Queensland, QLD 4072, 
Australia}


\begin{abstract}
The Knop--Sahi interpolation Macdonald polynomials are inhomogeneous and
nonsymmetric generalisations of the well-known Macdonald polynomials.
In this paper we apply the interpolation Macdonald polynomials to study
a new type of basic hypergeometric series of type $\gn$.
Our main results include a new $q$-binomial theorem, new $q$-Gauss sum, and 
several transformation formulae for $\gn$ series.
\end{abstract}

\maketitle

\tableofcontents

\newpage

\part{Interpolation Macdonald Polynomials}
\label{part1}

\section{Introduction}
\label{secIntro}

The Newton interpolation polynomials
\begin{equation}\label{NI}
N_k(x)=(x-x_1)\cdots(x-x_k)
\end{equation}
were used by Newton in his now famous expansion
\[
N(x)=\sum_{i=0}^k N(x_1)\partial_1\partial_2\cdots\partial_i N_i(x).
\]
Here $N(x)$ is an arbitrary polynomial of degree $k$ and 
$\partial_i$ (operators in this paper
act on the \emph{left}) is a Newton divided difference operator
\[
f(x_1,x_2,\dots)\partial_i=
\frac{f(\dots,x_{i+1},x_i,\dots)-f(\dots,x_i,x_{i+1},\dots)}{x_{i+1}-x_i}.
\]

Various multivariable generalisations of the Newton interpolation polynomials
exist in the literature, such as the Schubert polynomials \cite{Lascoux03}
and several types of Macdonald interpolation polynomials 
\cite{Knop97,Lascoux01,Okounkov97,Okounkov98,Okounkov98b,Sahi96,Sahi98}.
In this paper we are interested in the latter, providing generalisations
of \eqref{NI} when the interpolation points $x_1,\dots,x_k$ form a geometric 
progression
\[
(x_1,x_2,x_3,\dots)=(1,q,q^2,\dots).
\]
Then
\begin{equation}\label{NA}
N_k(x)=(x-1)(x-q)\cdots(x-q^{k-1}),
\end{equation}
and three equivalent characterisations may be given as follows.
\begin{enumerate}
\item $N_k(x)$ is the unique monic polynomial of 
degree $k$ such that $N_k(q^m)=0$ for $m\in\{0,1,\dots,k-1\}$.
\item $N_k(x)$ is the solution of the recurrence
\[
p_{k+1}(x)=q^k (x-1)p_k(x/q)
\]
with initial condition $p_0(x)=1$.
\item
Up to normalisation $N_k(x)$ is the unique polynomial eigenfunction, 
with eigenvalue $q^{-k}$, of the operator
\[
\xi=\tau\Bigl(1-\frac{1}{x}\Bigr)+\frac{1}{x},
\]
where $f(x)\tau=f(x/q)$.
\end{enumerate}

Knop \cite{Knop97} and Sahi \cite{Sahi96} 
generalised the Newton interpolation polynomials to a family 
of nonsymmetric, inhomogeneous polynomials
$M_u(x)$, labelled by compositions $u\in \N^n$ and depending
on $n$ variables; $x=(x_1,\dots,x_n)$.
The polynomials $M_u(x)$, known as the (nonsymmetric) Macdonald
interpolation polynomials or (nonsymmetric)
vanishing Macdonald polynomials,
form a distinguished basis in the ring $\Q(q,t)[x_1,\dots,x_n]$.
Remarkably, Knop and Sahi showed that all three characterisations 
of the Newton interpolation polynomials carry over to the multivariable
theory. What appears not to have been observed before, however,
is that the Macdonald interpolation polynomials may be employed
to build a multivariable theory of basic hypergeometric 
series of type $\gn$. For example, with $\MM_u(x)$ an
appropriate normalisation of $M_u(x)$, the following
$n$-dimensional extension of the famous $q$-binomial theorem
holds
\begin{equation}\label{qBT}
\sum_u a^{\abs{u}} \MM_u(x) =
\prod_{i=1}^n \frac{(at^{n-i})_{\infty}}{(ax_i)_{\infty}}.
\end{equation}

The $\gn$ basic hypergeometric series studied in the paper
are very different to existing multiple basic hypergeometric series, 
such as those pioneered by Gustafson and 
Milne \cite{Gustafson87,Milne85a,Milne85} 
or those studied subsequently by a large 
number of authors, see e.g., \cite{GR04,Milne01} and references therein.

\medskip

In Part~\ref{part1} of this paper, comprising of 
Sections~\ref{secIntro}--\ref{secE}, we lay the necessary 
groundwork for studying basic hypergeometric series based on the 
interpolation Macdonald polynomials. This in itself will involve the 
study of another type of multivariable basic hypergeometric series 
involving the function $\Skew{u}{v}{a}{b}$, which is a normalised 
connection coefficient between the interpolation polynomials 
$\MM_u(ax)$ and $\MM_v(bx)$.
An example of an identity for the connection coefficients is
the multivariable $q$-Pfaff--Saalsch\"utz sum
\[
\sum_v \frac{(a)_v}{(c)_v} \, \Skew{u}{v}{a}{b} \Skew{v}{w}{b}{c}
=\frac{(a)_w(b)_u}{(b)_w(c)_u} \, \Skew{u}{w}{a}{c}.
\]
Upon symmetrisation our identities for the function $\Skew{u}{v}{a}{b}$
generalise multiple series studied by the second author in \cite{Rains05}
in work on BC$_n$-symmetric polynomials.

In Part~\ref{part2}, containing Sections~\ref{secin}--\ref{secKTW}, 
we define the $\gn$ basic hypergeometric series and prove a number of 
important results, such as a multiple analogues of the $q$-binomial and 
$q$-Gauss sums for interpolation Macdonald polynomials.
By taking the top-homogeneous components of $\gn$ identities we also
obtain results for $\sln$ basic hypergeometric function involving the
nonsymmetric (homogeneous) Macdonald polynomials $E_u(x)$.

\section{Compositions}
\label{secCompositions}

Let $\N=\{0,1,2,\dots\}$.
In this paper the letters $u,v,w$ will be reserved for compositions
of length $n$, i.e., $u=(u_1,\dots,u_n)\in\N^n$.
Occasionally, when more than three compositions are required in one
equation, we also use $\bar{u},\bar{v},\bar{w}$, where $u$ and $\bar{u}$
are understood to be independent.
For brevity the trivial composition $(0,\dots,0)\in\N^n$ will be written
simply as $0$.
The sum of the parts of the composition $u$ is denoted by $\abs{u}$, i.e.,
$\abs{u}=u_1+\cdots+u_n$.

A composition is called \emph{dominant} if $u_i\geq u_{i+1}$ for all 
$1\leq i\leq n-1$, in other words, if $u$ is a \emph{partition}.
As is customary, we often use the Greek letters $\la,\mu$ and $\nu$
to denote partitions --- all of which are in $\N^n$ in this paper
(by attaching strings of zeros if necessary).

The symmetric group $\Symm_n$ acts on compositions by permuting the
parts. The unique partition in the $\Symm_n$-orbit of $u$ is denoted
by $u^{+}$. The staircase partition $\delta$ is defined as 
$\delta=(n-1,n-2,\dots,1,0)$, and $t^{\delta}$ is shorthand for 
$(t^{n-1},\dots,t,1)$. More generally we write 
$a^u b^v=(a^{u_1}b^{v_1},\dots,a^{u_n}b^{v_n})$
for $a,b$ scalars and $u,v$ compositions.

Given a composition $u$, we define its \emph{spectral vector} $\spec{u}$ by
\begin{equation}\label{spectral}
\spec{u}:=q^u t^{\delta\,\sigma_u},
\end{equation}
where $\sigma_u\in\Symm_n$ is the unique permutation of minimal length
such that $u=u^{+}\sigma_u$. Note that $\spec{\la}=q^{\la}t^{\delta}$
and $\spec{0}=t^{\delta}$.
Less formally, $\spec{u}$ is the unique permutation 
\[
(q^{u_1}t^{k_1},q^{u_2}t^{k_2},\dots,q^{u_n}t^{k_n})
\]
of $q^{u^{+}}t^{\delta}$ such that, if $u_i=u_j$ for $i<j$, then $k_i>k_j$.
For example, the spectral vector of $u=(2,4,2,0,1,2,1)$ is given by
\[
\spec{u}=(q^2t^5,q^4t^6,q^2t^4,1,qt^2,q^2t^3,qt).
\]

\medskip

The diagram of the composition $u$ is the set of points $(i,j)\in\Z^2$
such that $1\leq i\leq n$ and $1\leq j\leq u_i$. We write $v\subseteq u$ 
if the diagram of $v$ is contained in the diagram of $u$ (i.e., if
$v_i\leq u_i$ for all $1\leq i\leq n$) and $v\subset u$ if $v\subseteq u$
and $v\neq u$.
For $s=(i,j)$ in the diagram of $u$, the
arm-length $a(s)$, arm-colength $a'(s)$, leg-length $l(s)$
and leg-colength $l'(s)$ are given by \cite{Knop97,Sahi96b}
\[
a(s)=u_i-j, \qquad a'(s)=j-1
\]
and
\begin{align*}
l(s)&=\abs{\{k>i:~j\leq u_k\leq u_i\}}+\abs{\{k<i:~j\leq u_k+1\leq u_i\}} 
\\[1mm]
l'(s)&=\abs{\{k>i:~u_k>u_i\}}+\abs{\{k<i:~u_k\geq u_i\}} .
\end{align*}
For dominant compositions these last two definitions reduce to the
usual leg-length and leg-colength for partitions \cite{Macdonald95}.

Since $l'(s)$ only depends on the row-coordinate of $s$, we will
also write $l'(i)$. We may then use the above definition for
empty rows (i.e., rows such that $u_i=0$) as well. 
It is easily seen that
\begin{equation}\label{ulp}
\spec{u}_i=q^{u_i}t^{n-1-l'(i)}.
\end{equation}

Two statistics on compositions frequently used in this paper are
\[
n(u)=\sum_{s\in u} l(s)\And n'(u)=\sum_{s\in u} a(s)=
\sum_{i=1}^n \binom{u_i}{2}.
\]
Obviously, $n'(u)=n'(u^{+})$. Moreover, if $\la'$ denotes the conjugate 
of the partition $\la$, then $n'(\la)=n(\la')$.
The main reason for introducing $n'(u)$ is thus to avoid typesetting
$n((u^{+})')$ in many of our formulae.

\medskip

Throughout this paper $q\in\Complex$ is fixed such that $\abs{q}<1$.
Then 
\[
(b)_{\infty}=(b;q)_{\infty}:=\prod_{i=0}^{\infty} (1-bq^i)
\]
and
\[
(b)_k=(b;q)_k:=\frac{(b)_{\infty}}{(bq^k)_{\infty}}=
\prod_{i=0}^{k-1} (1-bq^i)
\]
are the standard $q$-shifted factorials \cite{GR04}.
The last equation is extended to compositions $u$ by
\begin{equation}\label{qsf}
(b)_u=(b;q,t)_u:=\prod_{s\in u}(1-b\, q^{a'(s)}t^{-l'(s)}).
\end{equation}
Note that this is invariant under permutations of $u$
\[
(b)_u=(b)_{u^{+}}=\prod_{i=1}^n (bt^{1-i};q)_{u^{+}_i}.
\]
Alternatively, by \eqref{ulp}, we can write
\begin{equation}\label{bu}
(bt^{n-1})_u=
\prod_{i=1}^n \frac{(b\spec{0}_i)_{\infty}}{(b\spec{u}_i)_{\infty}}.
\end{equation}

We also employ condensed notation for (generalised) $q$-shifted
factorials, setting
\[
(a_1,a_2,\dots,a_N)_u=\prod_{i=1}^N (a_i)_u.
\]

A special role in the theory of basic hypergeometric series 
is played by the $q$-shifted factorial $(q)_k$.
In the multivariable theory this role is played not by $(q)_u$,
but by
\[
c'_u=c'_u(q,t):=\prod_{s\in u} (1-q^{a(s)+1}t^{l(s)}).
\]
Occasionally we also need the related functions
\[
c_u=c_u(q,t):=\prod_{s\in u} (1-q^{a(s)}t^{l(s)+1})
\]
and
\begin{equation}\label{bdef}
b_u=b_u(q,t):=\frac{c_u}{c'_u}.
\end{equation}
For partitions these are standard in Macdonald polynomial theory, 
see \cite{Macdonald95}.

\section{Interpolation Macdonald polynomials}

Let $x=(x_1,\dots,x_n)$.
The definition of the interpolation Macdonald polynomial 
$M_u(x)=M_u(x;q,t)$ is deceptively simple 
\cite{Knop97,Lascoux01,Lascoux08,Sahi96,Sahi98}. 
It is the unique polynomial of degree $\abs{u}$ such that
\begin{equation}\label{Mudef}
M_u(\spec{v})=0 \quad\text{for $\abs{v}\leq \abs{u}$, $u\neq v$}
\end{equation}
and such that the coefficient of $x^u$ is $q^{-n'(u)}$.
Note that in the one-variable case
\[
M_u(x)=q^{-\binom{u}{2}} N_u(x),
\]
with $N_u(x)$ the Newton polynomial \eqref{NA}.
Comparing the above definition with those of Knop 
(denoted $E_u$ in \cite{Knop97}) and Sahi (denoted $G_u$ in \cite{Sahi98})
we find that
\[
M_u(x)=q^{-n'(u)}t^{(n-1)\abs{u}} F_u(xt^{1-n}), \qquad
F=E,G.
\]

\medskip

One of the key results in the theory is that the polynomials
$M_u$ can be computed recursively in much the same 
way as the Newton interpolation polynomials.
Let $s_i\in \Symm_n$ be the 
elementary transposition interchanging the 
variables $x_i$ and $x_{i+1}$.
Then the operator $T_i$ (acting on Laurent polynomials in $x$)
is defined as the unique operator that commutes with functions
symmetric in $x_i$ and $x_{i+1}$, such that
\[
1 T_i=t \And x_{i+1}T_i=x_i.
\]
More explicitly,
\[
T_i=t+(s_i-1)\,\frac{tx_{i+1}-x_i}{x_{i+1}-x_i}.
\]
It may readily be verified (see e.g., \cite{LS83,Lusztig89})
that the $T_i$ for $1\leq i\leq n-1$ satisfy the defining relations of
the Hecke algebra of the symmetric group
\begin{gather*}
T_iT_{i+1}T_i=T_{i+1}T_iT_{i+1}, \\[2mm]
T_iT_j=T_jT_i \quad\text{for $\abs{i-j}\neq 1$,} \\[2mm]
(T_i-t)(T_i+1)=0.
\end{gather*}
To describe the recursion we also require the operator $\tau$
\[
f(x)\tau:=f(x_n/q,x_1,\dots,x_{n-1})
\] 
and the raising operator $\phi$
\[
f(x)\phi:=f(x\tau)(x_n-1).
\]

According to \cite[Theorems 4.1 \& 4.2]{Knop97a} and
\cite[Corollary 4.4 \& Theorem 4.5]{Sahi96} 
the interpolation Macdonald polynomials can now be computed
as follows:
\[
M_0(x)=1,
\]
\[
M_{(u_2,\dots,u_{n-1},u_1+1)}(x)=M_u(x)\phi
\]
and 
\[
M_{us_i}(x)=M_u(x)\biggl(T_i+\frac{t-1}{\spec{u}_{i+1}/\spec{u}_i-1}\biggr)
\qquad \text{for $u_i<u_{i+1}$.}
\]
For a more general view of this recursive construction in terms of
Yang--Baxter graphs we refer the reader to \cite{Lascoux01}.

\medskip

A third description of the $M_u$ requires Knop's generalised
Cherednik operators $\Xi_i$ \cite{Knop97a,Knop97}
\begin{align*}
\Xi_i&:=t^{1-n} T_{i-1}\cdots T_1\tau(x_n-1)
T_{n-1}\cdots T_i \, \frac{1}{x_i}+\frac{1}{x_i} \\
&\hphantom{:}=
t^{1-i} T_{i-1}\cdots T_1\tau
\biggl(1-\frac{1}{x_n}\biggr)
T_{n-1}^{-1}\cdots T_i^{-1}+\frac{1}{x_i}
\end{align*}
for $1\leq i\leq n$. The $\Xi_i$ are mutually commuting and
the interpolation Macdonald polynomials
may be shown to be simultaneous eigenfunctions of the $\Xi_i$.
Specifically, \cite[Theorem 3.6]{Knop97}
\[
M_u(x) \Xi_i=M_u(x)/\spec{u}_i.
\]

Observe that the $T_i$ (and hence their inverses) and $\tau$
are degree preserving operators. The top-homogeneous degree $M^t(x)$ of
$M_u(x)$ thus satisfies
\[
M^t_u(x) Y_i^{-1}=M^t_u(x)/\spec{u}_i,
\]
where
\[
Y_i^{-1}=t^{1-i} T_{i-1}\cdots T_1\tau T_{n-1}^{-1}\cdots T_i^{-1}.
\]
Since the $Y_i$ are precisely the Cherednik operators \cite{Cherednik91},
which have the nonsymmetric Macdonald polynomials $E_u(x)=E(x;q,t)$ as 
simultaneous eigenfunctions with eigenvalues $\spec{u}_i$, it follows that
the top-homogeneous component of $M_u(x)$ is given by $E_u(x)$ 
\cite[Theorem 3.9]{Knop97}. To be more precise, since the coefficient of 
$x^u$ in $E_u(x)$ is $1$, it follows that
\begin{equation}\label{EfromM}
E_u(x)=q^{n'(u)} \lim_{a\to 0} a^{\abs{u}} M_u(x/a).
\end{equation}

\medskip

The $\gn$ basic hypergeometric series studied in Part~\ref{part2}
of this paper contain the Macdonald interpolation polynomials as 
key-ingredient.
In developing the theory we also frequently
need specific knowledge about $E_u(x)$.
Both these functions almost exclusively occur 
in combination with $c'_u$, and it will be convenient to
define the normalised Macdonald polynomials
\begin{equation}\label{MMdef}
\MM_u(x):=q^{n'(u)} t^{n(u)} \frac{M_u(x)}{c'_u}
\end{equation}
and
\begin{equation*}
\EE_u(x):=t^{n(u)} \frac{E_u(x)}{c'_u}.
\end{equation*}
Note that in the one-variable case
\begin{equation}\label{MMn1}
\MM_u(x)=x^u \, \frac{(1/x)_u}{(q)_u}\and
\EE_u(x)=\frac{x^u}{(q)_u},\quad u\in\N.
\end{equation}
Also, from \eqref{EfromM},
\begin{equation}\label{EEfromMM}
\EE_u(x)=\lim_{a\to 0} a^{\abs{u}} \MM_u(x/a).
\end{equation}
We use repeatedly (and implicitly) in subsequent sections that 
$\{\MM_u(x): \abs{u}\leq k\}$ (resp. $\{\EE_u(x): \abs{u}=k\}$)
forms a $\Q(q,t)$-basis in the space of polynomials
of degree $\leq k$ (exactly $k$) in $n$ variables.

On several occasions we also need the symmetric analogues of
$M_u$ and $E_u$, denoted $\MS_{\la}$ and $P_{\la}$, respectively.
$P_{\la}$ is of course the Macdonald polynomial as originally
introduced by Macdonald \cite{Macdonald88,Macdonald95}.
Using the same normalisation as before, i.e.,
\[
\MMS_{\la}(x)=t^{n(\la)} \frac{\MS_{\la}(x)}{c'_{\la}} \And
\PP_{\la}(x)=t^{n(\la)} \frac{P_{\la}(x)}{c'_{\la}},
\]
the symmetric polynomials are given by
\begin{equation}\label{symmetrise}
\MMS_{\la}(x):=\sum_{u^{+}=\la} \MM_u(x) \And \PP_{\la}(x):=
\sum_{u^{+}=\la} \EE_u(x).
\end{equation}
(This is of not the standard way to define the symmetric polynomials
but provides the most useful description for our purposes.)

A final result about the interpolation polynomials needed subsequently is
Sahi's principal specialisation formula \cite[Theorem 1.1]{Sahi98}.
Most convenient will be to normalise Sahi's formula by Cherednik's principal
specialisation formula for $\EE_u(x)$ \cite[Main Theorem]{Cherednik95}. Then
\begin{equation}\label{Mz0}
\MM_u(z\spec{0})=(1/z)_u \EE_u(z\spec{0}).
\end{equation}
For later reference we state the $z\to 0$ limit of this separately
\begin{equation}\label{M0E0}
\MM_u(0)=\tau_u \,\EE_u(\spec{0}),
\end{equation}
where $\MM_u(0)$ is shorthand for $\MM_u(0,\dots,0)$ and
$\tau_u=\tau_u(q,t)=\tau^{-1}_u(1/q,1/t)$ is defined as
\begin{equation}\label{tau}
\tau_u:=(-1)^{\abs{u}}q^{n'(u)}t^{-n(u^{+})}.
\end{equation}

\section{The Okounkov and Sahi binomial formulas}
\label{SecGqb}

Another important ingredient in the $\gn$ basic hypergeometric series
are Sahi's \cite{Sahi98} generalised $q$-binomial coefficients
\begin{equation}\label{genqbin}
\qbin{u}{v}=\qbin{u}{v}_{q,t}:=\frac{\MM_v(\spec{u})}{\MM_v(\spec{v})}.
\end{equation}
Note that, since $\MM_0(x)=1$, 
\[
\qbin{u}{0}=\qbin{u}{u}=1.
\]
{}From \eqref{Mudef} it follows that 
\begin{equation}\label{qbin0}
\qbin{u}{v}=0 \quad \text{if $\abs{u}\leq \abs{v}$, $u\neq v$.}
\end{equation}
A Theorem of Knop \cite[Theorem 4.5]{Knop97} also implies that
\begin{equation}\label{Knop0}
\qbin{u}{v}=0 \quad\text{if $v^{+}\not\subseteq u^{+}$.}
\end{equation}

Using \eqref{NA} it follows that in the one-variable case \eqref{genqbin}
reduces to the definition of the classical $q$-binomial coefficients
\begin{equation}\label{gausspoly}
\qbin{u}{v}=\frac{(q^{u-v+1})_v}{(q)_v}, \quad u,v\in\N,
\end{equation}
known to be polynomials in $q$ with positive integer coefficients.
For general $n$ and generic $u$ and $v$ the generalised $q$-binomial 
coefficients are, however, rational functions in $q$ and $t$. 
For example, if $u$ and 
\[
u^{(k)}:=(u_1,\dots,u_{k-1},u_k+1,u_{k+1},\dots,u_n)
\]
are both dominant, then \cite{Lascoux08}
\[
\qbin{u^{(k)}}{u}=
\frac{1-q^{u_k+1} t^{n-k}}{1-q}
\prod_{i=1}^{k-1} 
\frac{1-q^{u_k-u_i}t^{i-k-1}}{1-q^{u_k-u_i}t^{i-k}}
\prod_{i=k+1}^n 
\frac{1-q^{u_k-u_i+1}t^{i-k-1}}{1-q^{u_k-u_i+1}t^{i-k}}.
\]
Since \cite[page 14]{Lascoux08}
\[
\MM_u(\spec{u})=\tau_u\,t^{(n-1)\abs{u}}
\]
definition \eqref{genqbin} can also be written as
\begin{equation}\label{MMqbin}
\MM_v(\spec{u})=\tau_v \, t^{(n-1)\abs{v}}\qbin{u}{v}.
\end{equation}

\medskip

In \cite{Sahi98} Sahi proved a binomial formula for the
interpolation Macdonald polynomials $\MM_u(x)$ that will
be of importance later. In fact, we will be needing slightly more
than what may be found in Sahi's paper, and for clarity's sake
it is best to first recall the analogous results obtained by Okounkov
\cite{Okounkov97} for the symmetric interpolation 
polynomials $\MMS_{\la}(x)$. 

Before we can state Okounkov's theorem we need the symmetric
analogue of the generalised $q,t$-binomial coefficients 
\eqref{genqbin}, introduced independently by Lassalle and
Okounkov \cite{Lassalle98,Okounkov97} 
\[
\binom{\la}{\mu}=\binom{\la}{\mu}_{q,t}:=
\frac{\MMS_{\mu}(\spec{\la})}{\MMS_{\mu}(\spec{\mu})}.
\]
(Lassalle's definition is quite different
but may in fact be shown to be equivalent to the above.)
By symmetrisation it easily follows that for any composition
$u$ in the $\Symm_n$-orbit of $\la$ 
(i.e., any $u$ such that $u^{+}=\la$) \cite{BF00}
\[
\binom{\la}{\mu}=\sum_{v^{+}=\mu} \qbin{u}{v}.
\]
For later comparison with the nonsymmetric case it will also be useful 
to define 
\begin{equation}\label{MSp}
\MMS_{\la}'(x)=\MMS_{\la}'(x;q,t):=
\tau_{\la}^{-1} (q^{-1}t^{n-1})^{\abs{\la}}
\MMS_{\la}(t^{1-n}x;q^{-1},t^{-1}).
\end{equation}

\begin{theorem}[Okounkov's binomial formula]\label{thmO}
For $\la$ a partition of at most $n$ parts
\[
\MMS_{\la}(ax)=\sum_{\mu}
a^{\abs{\mu}}\binom{\la}{\mu}_{\!q^{-1},t^{-1}}
\frac{{\MMS_{\la}(a\spec{0})}}{\MMS_{\mu}(a\spec{0})}\,
\MMS_{\mu}'(x).
\]
\end{theorem}

Okounkov's theorem has a certain inversion symmetry as follows.
Upon replacing $(a,x,q,t)\mapsto (1/a,axt^{1-n},1/q,1/t)$
(note that $\spec{0}_{(q,t)\mapsto (1/q,1/t)}=t^{1-n}\spec{0}$)
and using \eqref{MSp} as well as \cite[page~537]{Okounkov97}
\[
\MMS_{\mu}(a^{-1}t^{1-n}\spec{0};q^{-1},t^{-1})=
(a^{-1}qt^{1-n})^{\abs{\mu}}\MMS_{\mu}(a\spec{0}),
\]
one finds
\begin{equation}\label{MSdual}
a^{\abs{\la}} \MMS_{\la}'(x)=
\sum_{\mu}  \frac{\tau_{\mu}}{\tau_{\la}}
\binom{\la}{\mu} \frac{\MMS_{\la}(a\spec{0})}{\MMS_{\mu}(a\spec{0})}\,
\MMS_{\mu}(ax).
\end{equation}
Replacing $(\la,\mu)\mapsto(\mu,\nu)$, multiplying by
$\binom{\la}{\mu}_{\! q^{-1},t^{-1}}\MMS_{\la}(a\spec{0})/
\MMS_{\mu}(a\spec{0})$ and then summing over $\mu$, results in
\[
\MMS_{\la}(ax)=
\sum_{\mu,\nu} \frac{\tau_{\nu}}{\tau_{\mu}}
\binom{\la}{\mu}_{\! q^{-1},t^{-1}} \binom{\mu}{\nu}
\frac{\MMS_{\la}(a\spec{0})}{\MMS_{\nu}(a\spec{0})}\,
\MMS_{\nu}(ax),
\]
where the sum over $\mu$ on the left has been performed using 
Okounkov's binomial formula. Equating coefficients of $\MMS_{\la}(ax)$ 
(and replacing $(q,t)\mapsto (1/q,1/t)$) results
in the orthogonality relation \cite[page 540]{Okounkov97}
\begin{equation}\label{deltasym}
\sum_{\mu}\frac{\tau_{\mu}}{\tau_{\la}} 
\binom{\la}{\mu}\binom{\mu}{\nu}_{\! q^{-1},t^{-1}}
=\delta_{\la\nu}.
\end{equation}

An alternative viewpoint is that Theorem~\ref{thmO} may be
inverted using \eqref{deltasym} but that, up to a change of
parameters, this inverted form is equivalent to the original
one. As we shall see shortly, this inversion symmetry is absent
in Sahi's nonsymmetric analogue of Theorem~\ref{thmO}.

Before we can state this result we need to introduce the nonsymmetric
analogue of the function $\MMS_{\la}'(x)$ as follows.
Extend definition \eqref{spectral} to all 
integral vectors $I\in\Z^n$ by $\spec{I}:=q^I t^{\delta\,\sigma_I}$,
where $\sigma_I\in\Symm_n$ is the unique permutation of minimal length
such that $I=I^{+}\sigma_u$. Here $I^{+}$ is a dominant integral vector,
i.e., $I^{+}_i\geq I^{+}_{i+1}$ for all $i$. Then
$\MM'_u(x)=\MM'_u(x;q,t)$ is the unique polynomial such that
the top-homogeneous terms (i.e., the terms of degree $\abs{u}$)
of $\MM'_u(x)$ and $\MM_u(x)$ coincide
(and are thus given by $\EE_u(x)$), and such that
\[
\MM'_u(\spec{-v})=0\qquad\text{for $\abs{v}<\abs{u}$.}
\]
We remark that, generally, $\MM'_u(\spec{-v})\neq 0$ if
$\abs{v}=\abs{u}$.
Unlike the symmetric case, $\MM'_u(x)$ does not simply follow
from $\MM_u(x)$ by the map $(q,t)\mapsto (1/q,1/t)$.

\begin{theorem}[Sahi's binomial theorem]\label{thmS1}
For $u\in\N^n$
\[
\MM_u(ax)=
\sum_v a^{\abs{v}}\qbin{u}{v}_{q^{-1},t^{-1}}
\frac{\MM_u(a\spec{0})}{\MM_v(a\spec{0})}\, \MM'_v(x).
\]
\end{theorem}
Replacing $x\mapsto x/a$ and letting $a$ tend to zero using
\[
\lim_{a\to 0} a^{\abs{u}} \MM'_u(x/a)=\EE_u(x)
\]
and \eqref{M0E0} expresses $\MM_u(x)$ in terms of $\EE_v(x)$.
\begin{corollary}
For $u\in\N^n$
\begin{equation}\label{ME}
\MM_u(x)=\sum_v \frac{\tau_u}{\tau_v}\,
\frac{\EE_u(\spec{0})}{\EE_v(\spec{0})}
\qbin{u}{v}_{q^{-1},t^{-1}} \EE_v(x).
\end{equation}
\end{corollary}

In the absence of a non-symmetric analogue of \eqref{MSp}
Sahi's result lacks the inversion symmetry of Theorem~\ref{thmO},
and by only minor changes of the proof given in \cite{Sahi98}
one can also obtain the nonsymmetric version of \eqref{MSdual}.
\begin{theorem}\label{thmS2}
For $u\in\N^n$ 
\[
a^{\abs{u}} \MM'_u(x)=
\sum_v \frac{\tau_v}{\tau_u} \qbin{u}{v}
\frac{\MM_u(a\spec{0})}{\MM_v(a\spec{0})} \, \MM_v(ax).
\]
\end{theorem}
The key to the proof is the following dual version of Sahi's
reciprocity theorem \cite[Theorem 1.2]{Sahi98}: there exists a unique 
polynomial $\RR_v(x)$ of degree not exceeding $\abs{v}$ 
with coefficients in $\Q(q,t,a)$, such that 
\[
\RR_v(\spec{u})=\tau_u \, a^{\abs{u}} 
\frac{\MM'_u(\spec{v}/a)}{\MM_u(a\spec{0})}
\]
for all $u\in\N^n$.
{}From this one can compute the coefficients $c_{uv}$ in 
\[
\tau_u \, a^{\abs{u}} \frac{\MM'_u(x/a)}{\MM_u(a\spec{0})}=
\sum_v c_{uv} \MM_v(x),
\]
see \cite{Sahi98} for more details.

Of course, the consistency of the last two theorems dictates that
\begin{equation}\label{inv}
\sum_v \frac{\tau_v}{\tau_u}\qbin{u}{v}\qbin{v}{w}_{q^{-1},t^{-1}}
=\delta_{uw}.
\end{equation}
Theorem~\ref{thmS2} and equation \eqref{inv} each imply the inverse 
of \eqref{ME}.
\begin{corollary}\label{CorS2}
For $u\in\N^n$
\begin{equation}\label{EM}
\EE_u(x)=
\sum_v \frac{\EE_u(\spec{0})}{\EE_v(\spec{0})}
\qbin{u}{v}\MM_v(x).
\end{equation}
\end{corollary}

\section{One-variable basic hypergeometric series}
\label{BHS}

Suppressing the $q$-dependence, the standard notation for single-variable 
basic hypergeometric series is~\cite{GR04}
\begin{equation}\label{g1series}
{_r\phi_s}\biggl[\genfrac{}{}{0pt}{}
{a_1,\dots,a_r}{b_1,\dots,b_s};x\biggr]:=
\sum_{k=0}^{\infty}
\frac{(a_1,\dots,a_r)_k}{(q,b_1,\dots,b_s)_k}
\Bigl((-1)^k q^{\binom{k}{2}} \Bigr)^{s-r+1} x^k.
\end{equation}
A series is called \emph{terminating} if only a finite number
of terms contribute to the sum, for example if $a_r=q^{-m}$
with $m$ a nonnegative integer.

Nearly all important summation and transformation formulas concern 
$_{r+1}\phi_r$ series. A $_{r+1}\phi_r$ series such that 
\begin{subequations}
\begin{align}\label{bal}
a_1\cdots a_{r+1}x&=b_1\cdots b_r, \\
x&=q  \label{xq}
\end{align}
\end{subequations}
is called \emph{balanced}. For reasons that will become clear in 
Part~\ref{part2} we somewhat relax this terminology and refer to a 
$_{r+1}\phi_r$ series as balanced if \eqref{bal} holds but 
not necessarily \eqref{xq}. (On page 70 of \cite{GR04} a series 
satisfying \eqref{bal} is referred to as \emph{a series of type II}, but we 
prefer a somewhat more descriptive adjective.)

In the next section and in Part~\ref{part2} we will generalise
most of the simple summation and transformation formulas for 
$_{r+1}\phi_r$ series, and for later reference we list a number
of one-variable series below. The reader may find proofs of all of the
identities in the book by Gasper and Rahman \cite{GR04}.

\medskip
\noindent\emph{$q$-Binomial theorem.}
The $q$-binomial theorem \cite[Equation (III.3)]{GR04} is one of the most
famous identities for basic hypergeometric series, discovered around 
1850 by a number of mathematicians including Cauchy and Heine,
\begin{equation}\label{qbt}
{_1\phi_0}\biggl[\genfrac{}{}{0pt}{}
{a}{\text{--}};x\biggr]=
\frac{(ax)_{\infty}}{(x)_{\infty}}, \qquad \abs{x}<1.
\end{equation}

\medskip
\noindent\emph{$q$-Gauss sum.}
The $q$-Gauss sum is a generalisation of the $q$-binomial theorem
due to Heine \cite[Equation (II.8)]{GR04}
\begin{equation}\label{qG}
{_2\phi_1}\biggl[\genfrac{}{}{0pt}{}
{a,b}{c};\frac{c}{ab}\biggr]=
\frac{(c/a,c/b)_{\infty}}{(c,c/ab)_{\infty}}, \qquad \abs{c/ab}<1.
\end{equation}
In our terminology, the series on the left is balanced.
To obtain the $q$-binomial theorem it suffices to replace
$b\mapsto c/ax$ and take the $c\to 0$ limit.

\medskip
\noindent\emph{$q$-Chu--Vandermonde sums.}
The first $q$-Chu--Vandermonde sum is simply the terminating
version of the $q$-Gauss sum \cite[Equation (II.7)]{GR04}
\begin{equation}\label{qcv}
{_2\phi_1}\biggl[\genfrac{}{}{0pt}{}{a,q^{-k}}{c};\frac{cq^k}{a}\biggr]=
\frac{(c/a)_k}{(c)_k}.
\end{equation}
By replacing $q\mapsto 1/q$ or by reversing the order of summation
a second $q$-Chu--Vandermonde sum follows \cite[Equation (II.6)]{GR04}
\begin{equation}\label{qcv2}
{_2\phi_1}\biggl[\genfrac{}{}{0pt}{}
{a,q^{-k}}{c};q\biggr]=\frac{(c/a)_k}{(c)_k}\, a^k.
\end{equation}
The binomial formulas of Okounkov and Sahi are easily seen to be 
generalisations of the $q$-Chu--Vandermonde sums. For example,
Theorem~\ref{thmS1} for $n=1$ is \eqref{qcv2} with $(a,c,k)\mapsto (x,1/a,u)$
and Theorem~\ref{thmS2} for $n=1$ is \eqref{qcv} with 
$(a,c,k)\mapsto (1/ax,1/a,u)$.

\medskip
\noindent\emph{$q$-Pfaff--Saalsch\"utz sum.}
At the top of the summations considered in this paper is the
$q$-Pfaff--Saalsch\"utz sum, first derived by Jackson
and generalising the $q$-Gauss sum \cite[Equation (II.12)]{GR04},
\begin{equation}\label{qps1}
{_3\phi_2}\biggl[\genfrac{}{}{0pt}{}
{a,b,q^{-k}}{c,abq^{1-k}/c};q\biggr]=\frac{(c/a,c/b)_k}{(c,c/ab)_k}.
\end{equation}
This sum is balanced in the traditional sense and yields the $q$-Gauss
in the large $k$ limit.

\medskip
\noindent\emph{Heine's transformations.}
There are three Heine transformations for ${_2\phi_1}$ series.
Of interest in this paper are only two of these.
The first one is \cite[Equation (III.2)]{GR04}
\begin{equation}\label{Heine1}
{_2\phi_1}\biggl[\genfrac{}{}{0pt}{}{a,b}{c};x\biggr]=
\frac{(c/a,ax)_{\infty}}{(c,x)_{\infty}}\,
{_2\phi_1}\biggl[\genfrac{}{}{0pt}{}{a,abx/c}{ax};\frac{c}{a}\biggr]
\end{equation}
for $\max\{\abs{x},\abs{c/a}\}<1$, and simplifies to the $q$-Gauss
sum \eqref{qG} when the balancing condition $x=c/ab$ is imposed.
The second one is Heine's $q$-analogue of Euler's $_2F_1$ transformation
\cite[Equation (III.3)]{GR04}
\begin{equation}\label{Heine2}
{_2\phi_1}\biggl[\genfrac{}{}{0pt}{}{a,b}{c};x\biggr]=
\frac{(abx/c)_{\infty}}{(c)_{\infty}}\,
{_2\phi_1}\biggl[\genfrac{}{}{0pt}{}{c/a,c/b}{c};\frac{abx}{c}\biggr]
\end{equation}
for $\max\{\abs{x},\abs{abx/c}\}<1$, and simplifies to the $q$-binomial
theorem \eqref{qbt} when $c=b$.

\medskip
\noindent\emph{$q$-Kummer--Thomae--Whipple formula.}
Sears' $q$-analogue of the Kummer--Thomae--Whipple formula
for $_2F_1$ series is given by \cite[Equation (III.9)]{GR04}
\begin{equation}\label{KTW1}
{_3\phi_2}\biggl[\genfrac{}{}{0pt}{}
{a,b,c}{d,e};\frac{f}{a}\biggr]=
\frac{(e/a,f)_{\infty}}{(f/a,e)_{\infty}}\,
{_3\phi_2}\biggl[\genfrac{}{}{0pt}{}
{a,d/b,d/c}{d,f};\frac{e}{a}\biggr],
\end{equation}
where $f=de/bc$ and $\max\{\abs{e/a},\abs{f/a}\}<1$.
Note that both $_3\phi_2$ series are balanced, and that
the limits $c,d\to\infty$ (such that $de/abc=x$) and
$a,e\to 0$ (such that $de/abc=x$) yield the 
Heine transformations \eqref{Heine1} and \eqref{Heine2} respectively.

\medskip

\noindent\emph{Sears' $_4\phi_3$ transformation.}
The most general result in our list is Sears' $_4\phi_3$
transformation formula~\cite[Equation (III.15)]{GR04}
\begin{equation}\label{Sears1}
{_4\phi_3}\biggl[\genfrac{}{}{0pt}{}{a,b,c,q^{-k}}{d,e,f};q\biggr]=
\frac{(e/a,f/a)_k}{(e,f)_k}\, a^n \,
{_4\phi_3}\biggl[\genfrac{}{}{0pt}{}{a,d/b,d/c,q^{-k}}
{d,aq^{1-n}/e,aq^{1-n}/f};q\biggr]
\end{equation}
for $abc=defq^{k-1}$. In the large $k$ limit this yields \eqref{KTW1} 
whereas for $d=c$ it simplifies to the $q$-Pfaff--Saalsch\"utz 
sum \eqref{qps1}.

\section{The function \texorpdfstring{$\Skew{u}{v}{a}{b}$}{}}
\label{secE}

A final ingredient needed in our study of $\gn$ basic hypergeometric series
are certain normalised connection coefficients between interpolation 
Macdonald polynomials.
\begin{definition}\label{defconnection}
For $u,v\in\N^n$ let the connection coefficients 
$c_{uv}(a,b)=c_{uv}(a,b;q,t)$ be given by 
\begin{equation}\label{cdef}
\MM_u(ax)=\sum_v c_{uv}(a,b) \MM_v(bx).
\end{equation}
Then
\begin{equation}\label{Skewdef}
\Skew{u}{v}{a}{b}=\Skewqt{u}{v}{a}{b}{q}{t}:=
\frac{\EE_v(\spec{0}/a)}{\EE_u(\spec{0}/b)}\, c_{uv}(a,b).
\end{equation}
\end{definition}

{}From the definition of the coefficients $c_{uv}(a,b)$
it immediately follows that 
\[
c_{uv}(a,b)=c_{uv}(ac,bc).
\]
Hence 
\begin{equation}\label{homEuv}
\Skew{u}{v}{ac}{bc}=c^{\abs{u}-\abs{v}} \Skew{u}{v}{a}{b}.
\end{equation}
Two other easy consequence of the definition are
\begin{equation}\label{deltaE}
\lim_{a\to b} \Skew{u}{v}{a}{b}=\delta_{uv}
\end{equation}
and the orthogonality relation
\begin{equation}\label{orthogonality}
\sum_v \Skew{u}{v}{a}{b}\Skew{v}{w}{b}{a}=\delta_{uw}.
\end{equation}

{}By \eqref{MMn1} and \eqref{qcv}
with $(a,c,k)\mapsto(1/bx,q^{1-u}a/b,u)$ it can be shown that
\begin{equation}\label{Euv1}
\Skew{u}{v}{1}{b}=(b)_{u-v}\qbin{u}{v} \quad \text{for $u,v\in\N$.}
\end{equation}
In general such a simple factorisation does not hold, although some 
of the features of \eqref{Euv1} lift to $n>1$.
\begin{lemma}\label{Lemqbin}
For $u,v\in\N^n$
\begin{subequations}
\begin{align}\label{Euvqbin}
\Skew{u}{v}{1}{0}&=\qbin{u}{v}, \\[2mm]
\Skew{u}{v}{0}{1}&=\frac{\tau_u}{\tau_v}\qbin{u}{v}_{q^{-1},t^{-1}} 
\label{Euvqbininv}
\end{align}
and
\begin{equation}\label{bu2}
\Skew{u}{0}{1}{b}=(b)_u.
\end{equation}
\end{subequations}
\end{lemma}
The connection coefficient $\Skew{u}{v}{a}{b}$ thus interpolates between 
generalised $q$-shifted factorials and
generalised $q$-binomial coefficients.

\begin{proof}
{}From \eqref{cdef} and \eqref{Skewdef} with $(a,b)\mapsto (1/a,1)$
\[
a^{\abs{u}}\MM_u(x/a)=\sum_v 
\frac{\EE_u(\spec{0})}{\EE_v(\spec{0})}\,
\Skew{u}{v}{1}{a} \MM_v(x),
\]
where we have also used \eqref{homEuv} to replace $\Skew{u}{v}{1/a}{1}$
by $a^{\abs{v}-\abs{u}} \Skew{u}{v}{1}{a}$.
By \eqref{EEfromMM} we can take the $a\to 0$ limit, giving
\[
\EE_u(x)=\sum_v \frac{\EE_u(\spec{0})}{\EE_v(\spec{0})}\,
\Skew{u}{v}{1}{0} \MM_v(x).
\]
Comparing this with \eqref{EM} completes the proof of \eqref{Euvqbin}.

{}From \eqref{cdef} and \eqref{Skewdef} with $(a,b)\mapsto (1,1/a)$
\[
\MM_u(x)=\sum_v a^{\abs{v}} \frac{\EE_u(\spec{0})}{\EE_v(\spec{0})}\,
\Skew{u}{v}{a}{1} \MM_v(x/a),
\]
where we have also used \eqref{homEuv} to replace $\Skew{u}{v}{1}{1/a}$
by $a^{\abs{v}-\abs{u}} \Skew{u}{v}{a}{1}$.
Taking the $a\to 0$ limit using \eqref{EEfromMM} yields
\[
\MM_u(x)=\sum_v \frac{\EE_u(\spec{0})}{\EE_v(\spec{0})}\,
\Skew{u}{v}{0}{1} \EE_v(x).
\]
Comparing this with \eqref{ME} results in \eqref{Euvqbininv}.

If we take $x=\spec{0}/ab$ in \eqref{cdef} and use the principal 
specialisation formula \eqref{Mz0} we get
\begin{equation}\label{bSkewa}
(b)_u=\sum_v \Skew{u}{v}{a}{b} (a)_v.
\end{equation}
Letting $a$ tend to $0$ using $\lim_{a\to 0} (a)_v=\delta_{v,0}$
results in \eqref{bu2}.
\end{proof}

The identity \eqref{bSkewa} is the $w=0$ instance of our next result.

\begin{proposition}[Multiple $q$-Chu--Vandermonde sum I]\label{PropqCV1}
For $u,w\in\N^n$
\begin{equation}\label{qCV1}
\sum_v \Skew{u}{v}{a}{b}\Skew{v}{w}{1}{a}=\Skew{u}{w}{1}{b}.
\end{equation}
\end{proposition}
For $n=1$ this is the $q$-Chu--Vandermonde sum \eqref{qcv} with 
$(c,k)$ replaced by $(q^{w-u+1}a/b,u-w)$.

Needed in our proof of the $\gn$ analogue of the $q$-Gauss sum
in Section~\ref{secGauss} is the following special case of \eqref{qCV1} 
obtained by substituting $b\to ab$, applying \eqref{homEuv} and 
letting $b\to 0$ using \ref{Euvqbin}.
\begin{corollary}\label{corcon}
For $u,w\in\N^n$
\[
\sum_v a^{\abs{u}-\abs{v}}\qbin{u}{v}\Skew{v}{w}{1}{a}=\qbin{u}{w}.
\]
\end{corollary}

\begin{proof}[Proof of Proposition \ref{PropqCV1}]
Double use of \eqref{cdef} gives
\[
\sum_v c_{uv}(a,b)c_{vw}(b,c)=c_{uw}(a,c).
\]
Hence, by \eqref{Skewdef},
\[
\sum_v \Bigl(\frac{a}{c}\Bigr)^{\abs{v}}
\Skew{u}{v}{a}{b}\Skew{v}{w}{b}{c}=
\Bigl(\frac{b}{c}\Bigr)^{\abs{u}}
\Bigl(\frac{a}{b}\Bigr)^{\abs{w}}
\Skew{u}{w}{a}{c}.
\]
Scaling $(a,b,c)\mapsto (ac,bc,abc)$ and applying \eqref{homEuv}
we obtain the desired result.
\end{proof}

\begin{proposition}\label{propER}
For $u,w\in\N^n$
\begin{subequations}
\begin{align}\label{eqER}
\Skew{u}{w}{a}{b}&=
\sum_v \frac{\tau_w}{\tau_v}\,\frac{(b)_u(a)_v}{(b)_v(a)_w}
\qbin{u}{v}_{q^{-1},t^{-1}} \qbin{v}{w} \\[1mm]
&=\sum_v a^{\abs{v}-\abs{w}}b^{\abs{u}-\abs{v}}
\frac{\tau_u}{\tau_v} 
\qbin{u}{v}_{q^{-1},t^{-1}} \qbin{v}{w}.
\label{eqERb}
\end{align}
\end{subequations}
\end{proposition}
These results not only give explicit expressions for the function
$\Skew{u}{w}{a}{b}$ but also show that it is polynomial
in $a$ and $b$, a fact that will be used later in our proof of 
Proposition~\ref{propDI}.

\begin{proof}[Proof of Proposition~\ref{propER}]
By \eqref{Mz0}, \eqref{M0E0} and \eqref{Skewdef},
and by substituting $(a,b)$ by $(1/b,1/a)$ using 
\[
\Skew{u}{v}{b^{-1}}{a^{-1}}=(ab)^{\abs{v}-\abs{u}} \Skew{u}{v}{a}{b},
\]
the first claim can be stated in equivalent form as
\[
c_{uw}(a,b)=
\sum_v \frac{\tau_w}{\tau_v}\,\Bigl(\frac{a}{b}\Bigr)^{\abs{v}}
\frac{\MM_u(a\spec{0})\MM_v(b\spec{0})}
{\MM_v(a\spec{0})\MM_w(b\spec{0})}
\qbin{u}{v}_{q^{-1},t^{-1}} \qbin{v}{w}.
\]
If we multiply this by $\MM_w(bx)$ and sum over $w$ then
\begin{align*}
\sum_w c_{uw}(a,b) \MM_w(bx)&=
\sum_{v,w} \frac{\tau_w}{\tau_v}\,\Bigl(\frac{a}{b}\Bigr)^{\abs{v}}
\frac{\MM_u(a\spec{0})\MM_v(b\spec{0})}
{\MM_v(a\spec{0})\MM_w(b\spec{0})}
\qbin{u}{v}_{q^{-1},t^{-1}} \qbin{v}{w} \MM_w(bx) \\
&=\sum_v a^{\abs{v}} \frac{\MM_u(a\spec{0})}
{\MM_v(a\spec{0})} \qbin{u}{v}_{q^{-1},t^{-1}} \MM'_v(x) \\[2mm]
&=\MM_u(ax).
\end{align*}
Here the second line follows from Theorem~\ref{thmS2}
and the third line from Theorem~\ref{thmS1}.

The second claim follows from the first by substituting
$(a,b)\mapsto(ac,bc)$ and then multiplying both sides by
$c^{\abs{w}-\abs{u}}$.
Using \eqref{homEuv} the large $c$ limit can now be taken
leading to \eqref{eqERb}.
\end{proof}

Our next result generalises the $q$-Pfaff--Saalsch\"utz sum \eqref{qps1},
and will be applied in Section~\ref{secKTW} to prove a $\gn$ version of 
the $q$-Kummer--Thomae--Whipple formula \eqref{KTW1}.

\begin{theorem}[Multiple $q$-Pfaff--Saalsch\"utz sum]\label{thmER2}
For $u,w\in\N^n$
\begin{equation}\label{qPS}
\sum_v \frac{(a)_v}{(c)_v}\,\Skew{u}{v}{a}{b}\Skew{v}{w}{b}{c}
=\frac{(a)_w(b)_u}{(b)_w(c)_u}\,\Skew{u}{w}{a}{c}.
\end{equation}
\end{theorem}
When $n=1$ this is \eqref{qps1} with $(a,b,c,k)\mapsto (aq^w,c/b,cq^w,u-w)$,
and for $b\to\infty$ this is (recall \eqref{Euvqbininv}) \eqref{eqER} 
with $b\mapsto c$.

\begin{proof}
Twice using \eqref{eqER} we find that
\begin{multline*}
\sum_v \frac{(a)_v(c)_w}{(b)_u(d)_v}\,
\Skew{u}{v}{a}{b} \Skew{v}{w}{c}{d} \\
=\sum_{v,\bar{v},\bar{w}} 
\frac{\tau_v\tau_w}{\tau_{\bar{v}}\tau_{\bar{w}}} \,
\frac{(a)_{\bar{v}}}{(b)_{\bar{v}}}\,
\frac{(c)_{\bar{w}}}{(d)_{\bar{w}}}
\qbin{u}{\bar{v}}_{q^{-1},t^{-1}} \qbin{\bar{v}}{v} 
\qbin{v}{\bar{w}}_{q^{-1},t^{-1}} \qbin{\bar{w}}{w}.
\end{multline*}
Summing over $v$ using \eqref{inv}, then performing the
trivial sum over $\bar{w}$, and finally replacing $\bar{v}$ by $v$
leads to
\[
\sum_v \frac{(a)_v(c)_w}{(b)_u(d)_v}\,
\Skew{u}{v}{a}{b} \Skew{v}{w}{c}{d} \\
=\sum_v
\frac{\tau_w}{\tau_v} \, \frac{(a,c)_v}{(b,d)_v}
\qbin{u}{v}_{q^{-1},t^{-1}} \qbin{v}{w}.
\]
The right-hand side is invariant under the maps $a\leftrightarrow c$ 
or $b\leftrightarrow d$. Hence,
\begin{subequations}
\begin{align}\label{abcd}
\sum_v \frac{(a)_v(c)_w}{(b)_u(d)_v}\,
\Skew{u}{v}{a}{b} \Skew{v}{w}{c}{d}
&=\sum_v \frac{(c)_v(a)_w}{(b)_u(d)_v}\,
\Skew{u}{v}{c}{b} \Skew{v}{w}{a}{d} \\
&=\sum_v \frac{(c)_v(a)_w}{(d)_u(b)_v}\,
\Skew{u}{v}{c}{d} \Skew{v}{w}{a}{b}.
\label{abcd2}
\end{align}
\end{subequations}
Setting $b=c$ in \eqref{abcd} using \eqref{deltaE} completes the proof.
\end{proof}

Several further multiple $q$-Chu--Vandermonde sums follow as limiting 
cases of Theorem~\ref{thmER2}.
\begin{corollary}[Multiple $q$-Chu--Vandermonde sums II--IV]\label{corqcv}
For $u,w\in\N^n$
\begin{subequations}
\begin{equation}\label{qCV2}
\sum_v \Skew{u}{v}{1}{a}\Skew{v}{w}{a}{b}=\Skew{u}{w}{1}{b}
\end{equation}
\begin{equation}\label{qCV3}
\sum_v \frac{\tau_v}{\tau_w} \, \frac{(b)_u}{(b)_v}
\qbin{u}{v} \Skew{v}{w}{a}{b}=\frac{(a)_u}{(a)_w} \qbin{u}{w}
\end{equation}
and
\begin{equation}\label{qCV4}
\sum_v b^{\abs{v}-\abs{w}} \frac{(a)_v}{(a)_w}\,
\Skew{u}{v}{a}{b}
\qbin{v}{w} =a^{\abs{u}-\abs{w}} \frac{(b)_u}{(b)_w} \qbin{u}{w}.
\end{equation}
\end{subequations}
\end{corollary}
For $n=1$ \eqref{qCV2} reduces to \eqref{qcv2} with 
$(a,c,k)\mapsto (b/a,q^{w-u+1}/a,u-w)$,
\eqref{qCV3} to \eqref{qcv} with
$(a,c,k)\mapsto (aq^w,q^{w-u+1}a/b,u-w)$
and \eqref{qCV4} to \eqref{qcv2} with 
$(a,c,k)\mapsto (b/a,bq^w,u-w)$.
We also remark that if we replace $a\mapsto 1/b$ in Theorem~\ref{thmS2},
then specialise $x=a\spec{0}$ using \eqref{Mz0} and \eqref{bu2}, we obtain
\[
\sum_v \frac{\tau_v}{\tau_u}\,
\frac{(b)_u}{(b)_v}\,\Skew{u}{0}{a}{b} \qbin{u}{v}
=\frac{\MM'_u(a\spec{0})}{\EE_u(\spec{0})}.
\]
Comparing this with the $w=0$ case of \eqref{qCV3} yields the
principal specialisation formula
\[
\MM'_u(a\spec{0})=\tau_u^{-1} (a)_u \EE_u(\spec{0}).
\]

\begin{proof}[Proof of Corollary~\ref{corqcv}]
Equation \eqref{qCV2} follows by first rescaling the parameters in
\eqref{qPS} as $(a,b,c)\mapsto (c,ac,bc)$ and then taking the $c\to 0$ limit.
Equation \eqref{qCV3} is \eqref{qPS} in the limit $(a,b,c)\to (\infty,a,b)$,
whereas \eqref{qCV4} corresponds to the $c\to 0$ limit of \eqref{qPS}.
\end{proof}

\begin{proposition}[Multiple Sears transformation]\label{propSears}
For $u,w\in\N^n$
\begin{multline}\label{Searsn}
\sum_v \frac{(aq/b,aq/c)_u(d,e)_v}{(aq/b,aq/c)_v(d,e)_w} \,
\Skew{u}{v}{1}{aq/de}\Skew{v}{w}{aq/de}{a^2q^2/bcde} \\
= \sum_v \frac{(aq/b,aq/d)_u(c,e)_v}{(aq/b,aq/d)_v(c,e)_w} \,
\Skew{u}{v}{1}{aq/ce} \Skew{v}{w}{aq/ce}{a^2q^2/bcde}.
\end{multline}
\end{proposition}
When $n=1$ this is \eqref{Sears1} with
\[
(a,b,c,d,e,f,k)\mapsto(eq^w,dq^w,aq/bc,aq^{w+1}/b,aq^{w+1}/c,q^{w-u}de/a,u-w).
\]
Later we shall also encounter the limiting case of \eqref{Searsn} obtained
by taking the $e\to\infty$ limit and relabelling the remaining variables
\begin{equation}\label{Searsn2}
\sum_v \tau_v\, \frac{(d,e)_u(b)_v}{(d,e)_v(b)_w}\qbin{u}{v} 
\Skew{v}{w}{a}{ac}
=\sum_v \tau_v\, \frac{(d,a)_u(d/c)_v}{(d,a)_v(d/c)_w}\qbin{u}{v}
\Skew{v}{w}{e}{ac},
\end{equation}
where $abc=de$.

\begin{proof}[Proof of Proposition~\ref{propSears}]
First we observe that thanks to \eqref{homEuv} the $q$-Pfaff--Saal\-sch\"utz
sum \eqref{qPS} can also be written with an additional parameter $d$ as
\begin{equation}\label{qps}
\sum_v \frac{(a)_v}{(c)_v}\,\Skew{u}{v}{ad}{bd}\Skew{v}{w}{bd}{cd}
=\frac{(a)_w(b)_u}{(b)_w(c)_u}\,\Skew{u}{w}{ad}{cd}.
\end{equation}
Making the substitutions $(a,b,c,d,u,v)\mapsto (c,d,aq/b,aq/cde,v,\bar{v})$
this leads to
\begin{multline*}
\Skew{v}{w}{aq/de}{a^2q^2/bcde} \\
=\frac{(d)_w(aq/b)_v}{(c)_w(d)_v}
\sum_{\bar{v}} \frac{(c)_{\bar{v}}}{(aq/b)_{\bar{v}}} \,
\Skew{v}{\bar{v}}{aq/de}{aq/ce}
\Skew{\bar{v}}{w}{aq/ce}{a^2q^2/bcde}.
\end{multline*}

Inserting the above expansion and interchanging the order of the $v$ 
and $\bar{v}$ sums we find
\begin{multline*}
\text{LHS}\eqref{Searsn}=
\sum_{\bar{v}}
\frac{(aq/b,aq/c)_u}{(c,e)_w} \,
\frac{(c)_{\bar{v}}}{(aq/b)_{\bar{v}}} \,
\Skew{\bar{v}}{w}{aq/ce}{a^2q^2/bcde} \\
\times \sum_v \frac{(e)_v}{(aq/c)_v} \,
\Skew{u}{v}{1}{aq/de} \Skew{v}{\bar{v}}{aq/de}{aq/ce}.
\end{multline*}
The sum over $v$ can now be performed by \eqref{qps} with
\[
(a,b,c,d)\mapsto (e,aq/d,aq/c,1/e)
\]
resulting in the right-hand side of the multiple Sears transform.
\end{proof}

\begin{proposition}[Duality, type I]\label{propDI}
For $k\in\N$ and $u,w\in\N^n$
\[
\sum_{\abs{v}=k}
\Skew{u}{v}{a}{b}\Skew{v}{w}{c}{d}=
\sum_{\abs{v}=\abs{u}+\abs{w}-k}
\Skew{u}{v}{c}{d}\Skew{v}{w}{a}{b}
\]
and
\[
\sum_{\abs{v}=k} \qbin{u}{v}\qbin{v}{w}
=\sum_{\abs{v}=\abs{u}+\abs{w}-k} \qbin{u}{v}\qbin{v}{w}.
\]
\end{proposition}

\begin{proof}
Replace $(a,b,c,d)\to (ae,be,ce,de)$ in \eqref{abcd2}, use
\eqref{homEuv}, and let $e$ tend to $0$. Then
\[
\sum_v \Skew{u}{v}{a}{b} \Skew{v}{w}{c}{d}
=\sum_v \Skew{u}{v}{c}{d} \Skew{v}{w}{a}{b}.
\]
From \eqref{eqERb} it follows that $\Skew{u}{v}{a}{b}$,
when viewed as a polynomial in $a$ and $b$, is homogeneous of 
degree $\abs{u}-\abs{v}$. If we therefore read off the term of degree
$k+\abs{w}$ in $c$ and $d$ in the above transformation
the first claim follows.

The second identity follows by taking $b=d=0$ and using \eqref{homEuv}
and \eqref{Euvqbin}.
\end{proof}

\begin{proposition}[Duality, type II]
For $u,v\in\N^n$
\begin{equation}\label{Skewdual}
\Skew{u}{v}{a}{b}=\frac{\tau_u}{\tau_v}\,
\Skewqt{u}{v}{b}{a}{q^{-1}}{t^{-1}}.
\end{equation}
\end{proposition}
Together with \eqref{orthogonality} this implies a generalisation 
of \eqref{inv}.

\begin{proof}
If we take the $b\to 0$ limit in \eqref{qPS} using \eqref{homEuv}
and Lemma~\ref{Lemqbin}, and replace $c\mapsto b$, we obtain
\[
\Skew{u}{w}{a}{b}
=\frac{a^{\abs{u}}}{b^{\abs{w}}} \, \frac{(b)_u}{(a)_w}
\sum_v \frac{\tau_v}{\tau_w} \Bigl(\frac{b}{a}\Bigr)^{\abs{v}}
\frac{(a)_v}{(b)_v}\qbin{u}{v} \qbin{v}{w}_{q^{-1},t^{-1}}.
\]
If we compare this with \eqref{eqER}, use \eqref{homEuv} and
\[
(a)_u=\tau_u\,a^{\abs{u}} (a^{-1};q^{-1},t^{-1})_u,
\]
the claim follows.
\end{proof}

Let $\g^u_{vw}=\g^u_{vw}(q,t)$ be the structure constants of the 
normalised nonsymmetric Macdonald polynomials
\begin{equation}\label{struc}
\EE_v(x)\EE_w(x)=\sum_u \g^u_{vw}\EE_u(x).
\end{equation}

\begin{proposition}\label{propEg}
For $u,v\in\N^n$
\begin{equation}\label{Eg}
\Skew{u}{v}{a}{b}=a^{\abs{u}-\abs{v}}\sum_w (b/a)_w\, \g^u_{vw}. 
\end{equation}
\end{proposition}
For $(a,v)=(1,0)$ this is \eqref{bu2} since $\g^u_{0,w}=\delta_{uw}$.
With some considerable pain \eqref{Eg} can be proved using 
the $q$-Pfaff--Saalsch\"utz sum \eqref{qPS}, but since it follows 
as an easy corollary of Theorem~\ref{thmGauss} we omit a proof here.

\medskip

We conclude this section with several remarks about the
function $\Skew{u}{v}{a}{b}$ and its symmetric counterpart.
For this purpose we first introduce $\la$-ring notation \cite{Lascoux03}, 
which also plays a crucial role in our proof of the $\gn$ $q$-binomial 
theorem in the next section.

For $\AA$ an alphabet (i.e., countable set with elements referred to
as letters) let $\abs{\AA}$ be its cardinality. 
If we adopt the usual additive notation for alphabets, that is,
$\AA=\sum_{a\in \AA} a$, then the disjoint union  and Cartesian product
of $\AA$ and $\BB$ may be written as
\[
\AA+\BB=\sum_{a\in\AA}a+\sum_{b\in\BB}b \and 
\AA\BB=\sum_{a\in\AA}\sum_{b\in\BB} ab.
\]
For $f$ a symmetric function, we define
\[
f[\AA]=f(a_1,a_2,\dots)\quad \text{for $\AA=a_1+a_2+\cdots$,}
\]
where $[\,$,~$]$ are referred to as plethystic brackets.
Let
\[
\sigma_z[\AA]=\prod_{a\in\AA}\frac{1}{1-za}.
\]
Then the complete symmetric function $S_k[\AA]$ is defined by its 
generating function as
\begin{equation}\label{complete}
\sigma_z[\AA]=\sum_{k=0}^{\infty} z^k S_k[\AA].
\end{equation}
More generally we define the complete symmetric function 
(and thus any symmetric function)
of the difference of two alphabets as
\begin{equation}\label{AminB}
\sigma_z[\AA-\BB]=\frac{\prod_{b\in\BB}(1-zb)}{\prod_{a\in\AA}(1-za)}.
\end{equation}
Hence, if $1/(1-q)$ denotes the infinite alphabet
$1+q+q^2+\cdots$, 
\begin{equation}\label{sigmaz}
\sigma_z\biggl[\frac{\AA-\BB}{1-q}\biggr]=
\frac{\prod_{b\in\BB}(zb)_{\infty}}{\prod_{a\in\AA}(za)_{\infty}}.
\end{equation}

Below we need \eqref{sigmaz} with $q\mapsto t$ and $\AA=a$, $\BB=b$ 
single-letter alphabets. More specifically, we consider 
the function
\[
P_{\la/\mu}\biggl[\frac{a-b}{1-t}\biggr],
\]
where $P_{\la/\mu}(x)=P_{\la/\mu}(x;q,t)$ is a skew Macdonald 
polynomial defined by
\[
P_{\la}[\AA+\BB]=\sum_{\mu} P_{\la/\mu}[\AA]P_{\mu}[\BB].
\]
If $f^{\la}_{\mu\nu}=f^{\la}_{\mu\nu}(q,t)$ are the 
$q,t$-Littlewood--Richardson coefficients
\[
P_{\mu}(x)P_{\nu}(x)=\sum_{\la} f^{\la}_{\mu\nu} P_{\la}(x),
\]
then, by \cite[page 344]{Macdonald95},
\[
P_{\la/\mu}(x)=
\sum_{\nu} \frac{b_{\mu}b_{\nu}}{b_{\la}}\, f^{\la}_{\mu\nu} P_{\nu}(x),
\]
where $b_{\la}$ is defined in \eqref{bdef}.
Combining this with \cite[page 338]{Macdonald95}
\[
P_{\la}\biggl[\frac{1-a}{1-t}\biggr]=t^{n(\la)} 
\frac{(a)_{\la}}{c_{\la}},
\]
it follows that
\begin{equation}\label{Plamu}
P_{\la/\mu}\biggl[\frac{a-b}{1-t}\biggr]
=a^{\abs{\la}-\abs{\mu}}\frac{b_{\mu}}{b_{\la}}
\sum_{\nu} t^{n(\nu)}
\frac{(b/a)_{\nu}}{c'_{\nu}} f_{\mu\nu}^{\la},
\end{equation}
where we have also used the homogeneity of $P_{\la/\mu}(x)$ and
the fact that $f_{\mu\nu}^{\la}=0$ unless $\abs{\la}=\abs{\mu}+\abs{\nu}$.
If we now define
\[
\f^{\la}_{\mu\nu}:=
t^{n(\mu)+n(\nu)-n(\la)} \frac{c'_{\la}}{c'_{\mu}c'_{\nu}}\, 
f^{\la}_{\mu\nu},
\]
so that (compare with \eqref{struc})
\[
\PP_{\mu}(x)\PP_{\nu}(x)=\sum_{\la} \f^{\la}_{\mu\nu}\PP_{\la}(x),
\]
and
\[
\PP_{\la/\mu}\biggl[\frac{a-b}{1-t}\biggr]
:=t^{n(\mu)-n(\la)} \frac{c_{\la}}{c_{\mu}}\,
P_{\la/\mu}\biggl[\frac{a-b}{1-t}\biggr],
\]
then \eqref{Plamu} simplifies to
\begin{equation}\label{PPf}
\PP_{\la/\mu}\biggl[\frac{a-b}{1-t}\biggr]
=a^{\abs{\la}-\abs{\mu}}
\sum_{\nu} (b/a)_{\nu}\, \f_{\mu\nu}^{\la}.
\end{equation}
(The reader is warned that the above choice of normalisations implies that
\[
\PP_{\la/0}[(a-b)/(1-t)]=t^{-2n(\la)} c_{\la}c'_{\la} \,\PP_{\la}[(a-b)/(1-t)]
\]
and not $\PP_{\la/0}[~]=\PP_{\la}[~]$.)
Comparing \eqref{PPf} with \eqref{Eg} we are led to conclude that
$\Skew{u}{v}{a}{b}$ is the nonsymmetric analogue of 
$\PP_{\la/\mu}[(a-b)/(1-t)]$. Indeed, from
\eqref{symmetrise} it follows that
\[
\PP_{\la/\mu}\biggl[\frac{a-b}{1-t}\biggr]
=\sum_{v^{+}=\mu}\Skew{u}{v}{a}{b},
\]
where $u^{+}=\la$. This in turn implies that the multiple $q$-Saalsch\"utz
sum \eqref{qPS} and Sears transformation \eqref{Searsn} are
the nonsymmetric analogues of \cite[Corollaries 4.9 \& 4.8]{Rains05}
respectively.

Up to trivial factors the generalised $q$-binomial coefficients arise as 
evaluations of the interpolation Macdonald polynomials $\MM_u(x)$, 
see \eqref{genqbin} or \eqref{MMqbin}.
It is thus natural to ask for generalisations of the $\MM_u(x)$ that yield
the $\Skew{u}{v}{a}{b}$ upon evaluation. In fact, in the symmetric theory 
such functions have already been constructed in \cite{Rains06}.
Specifically, $\Symm_n$-invariant \textit{rational functions} 
$f_{\mu}(x;a,b;q,t)=f_{\mu}(x;a,b)$ were defined
satisfying the following three conditions. (In comparison with \cite{Rains06}
the roles of $a$ and $b$ have been interchanged.) (1) The function
\[
\prod_{i=1}^n (qt^{n-1}a/x_i)_{\mu_1} f_{\mu}(x;a,b)
\]
is holomorphic in $\Complex^{\ast}$, (2) 
\[
f_{\mu}(b\spec{\la};a,b)=0 \quad\text{unless $\mu\subseteq\la$},
\]
and (3)
\[ 
\lim_{x\to a\spec{\la}} \prod_{i=1}^n (qt^{n-1}a/x_i)_{\la_1}
f_{\mu}(x;a,b)=0 \quad\text{unless $\la\subseteq\mu$}.
\]
Up to normalisation this uniquely fixes the functions $f_{\mu}(x;a,b)$.
Moreover,
\[
\PP_{\la/\mu}\biggl[\frac{a-b}{1-t}\biggr]=
f_{\mu}(b\spec{\la};a,b)\, d_{\la} e_{\mu},
\]
where $d_{\la}$ and $e_{\mu}$ are simple factors.
In the nonsymmetric theory similar rational interpolation functions
can be defined, which generalise the $\MM_u$ and which yield
the $\Skew{u}{v}{a}{b}$ upon evaluation. We hope to report on these
functions in a future publication.

\part{\texorpdfstring{$\gn$}{} Basic Hypergeometric Series}\label{part2}
Finally everything is in place to develop the theory of $\gn$ basic
hypergeometric series based on interpolation Macdonald polynomials. 

\section{Introduction}
\label{secin}

The multiple-series identities of the previous section
involving the connection coefficients $\Skew{u}{v}{a}{b}$
generalise all of the terminating identities listed in
Section~\ref{BHS}.
The identities for $\gn$ basic hypergeometric series proved
in the next few sections will generalise the non-terminating
identities in the list, except for \eqref{Heine2}.

Below two different types of $\gn$ basic hypergeometric series will 
be considered. At the top level are series of the form
\begin{equation}\label{gnseries}
\sum_u \frac{(a_1,\dots,a_{r-1})_u}{(b_1,\dots,b_{r-1},b_rt^{n-1})_u}\,
\Skew{u}{v}{c}{d}\MM_u(x),
\end{equation}
where $u,v\in\N^n$, and where the parameters satisfy 
the \emph{$\gn$ balancing condition}
\begin{equation}\label{gnbalanced}
a_1\cdots a_{r-1} d=b_1\cdots b_r.
\end{equation}
{}From \eqref{MMn1} and \eqref{Euv1} it follows that \eqref{gnseries}
for $n=1$ simplifies to
\[
x^v \frac{(a_1,\dots,a_{r-1},1/x)_v}{(q,b_1,\dots,b_r)_v} \,
{_{r+1}\phi_r}\biggl[\genfrac{}{}{0pt}{}
{a_1q^v,\dots,a_{r-1}q^v,q^v/x,d/c}{b_1q^v,\dots,b_rq^v};cx\biggr]
\]
which for $v=0$ is simply
\[
{_{r+1}\phi_r}\biggl[\genfrac{}{}{0pt}{}
{a_1,\dots,a_{r-1},1/x,d/c}{b_1,\dots,b_r};cx\biggr].
\]
Note that the $\gn$ balancing condition \eqref{gnbalanced} implies that the
$_{r+1}\phi_r$ series are balanced in the sense of \eqref{bal}.

At a lower level are $\gn$ series of the type
\begin{equation}\label{gnseriesa}
\sum_u c^{\abs{u}}\frac{(a_1,\dots,a_{r-1})_u}{(b_1,\dots,b_{r-1})_u}
\qbin{u}{v} \MM_u(x),
\end{equation}
obtained from \eqref{gnseries} by taking $d=b_r=0$,
so that no balancing condition holds.

In this paper we will consider the series \eqref{gnseries} and 
\eqref{gnseriesa} (and some related series discussed below) as formal
power series. However, for $\max\{\abs{cx_1},\dots,\abs{cx_n}\}$ 
sufficiently small the $\gn$ series may also be viewed as $\gn$ functions.

Typically, if a classical summation or transformation formula for
nonterminating basic hypergeometric series is balanced
(like the $q$-Gauss sum or the $q$-Kummer--Thomae--Whipple transformation)
it admits a $\gn$ generalisation involving series of the type
\eqref{gnseries}.
If, however, the parameters in the identity can be chosen freely (like in
the $q$-binomial theorem or Heine's $_2\phi_1$ transformations) it, at best,
admits a $\gn$ generalisation involving series of the type \eqref{gnseriesa}.

\medskip

Replacing $(b_r,c,d,x)\mapsto (b_re,ce,de,x/e)$ in \eqref{gnseries}, 
and using \eqref{EEfromMM} and \eqref{homEuv},
it follows that
\begin{multline*}
\lim_{c\to 0}
e^{\abs{v}}\sum_u\frac{(a_1,\dots,a_{r-1})_u}
{(b_1,\dots,b_{r-1},b_ret^{n-1})_u}\,
\Skew{u}{v}{ec}{ed}\MM_u(x/e) \\
=\sum_u \frac{(a_1,\dots,a_{r-1})_u}
{(b_1,\dots,b_{r-1})_u}\,\Skew{u}{v}{c}{d}\EE_u(x).
\end{multline*}
We will refer to the series on the right (or special cases thereof)
as $\sln$ basic hypergeometric series. 
Obviously, by the above limiting procedure, any $\gn$ series identity 
implies a corresponding $\sln$ series identity.

\medskip

Replacing $(b_r,c,d,u,v)\mapsto (b_rt^{1-n},ct^{1-n},dt^{1-n},v,w)$
in \eqref{gnseries}, multiplying by $t^{(n-1)\abs{w}}$
and specialising $x=\spec{u}$ we obtain the terminating series
\[
\sum_v \tau_v\, 
\frac{(a_1,\dots,a_{r-1})_v}{(b_1,\dots,b_r)_v}\,
\qbin{u}{v} \Skew{v}{w}{c}{d}
\]
subject to the balancing condition \eqref{gnbalanced}.
The $q$-Chu--Vandermonde sum \eqref{qCV3} and the transformation
\eqref{Searsn2} provide example of the above series.
We shall find in the next few sections that both \eqref{qCV3}
and \eqref{Searsn2} arise by the above type of specialisation
from $\gn$ basic hypergeometric series identities.

\section{A \texorpdfstring{$\gn$ $q$}{}-binomial theorem}
\label{SecAL}

Our first main result is a generalisation of the celebrated $q$-binomial
theorem \eqref{qbt}.
An equivalent way of stating this theorem follows from the substitutions 
$(a,x)\mapsto (q^v/x,ax)$, where $v$ is a nonnegative
integer. Then, after a shift in the summation index,
\[
\sum_{u=0}^{\infty}(ax)^u \qbin{u}{v} 
\frac{(1/x)_u}{(q)_u}=
(ax)^v \frac{(1/x)_v}{(q)_v}\,\frac{(aq^v)_{\infty}}{(ax)_{\infty}},
\]
where $\qbin{u}{v}$ is the $q$-binomial coefficient \eqref{gausspoly}.
Recalling \eqref{MMn1}, this may also be stated as
\[
\sum_{u=0}^{\infty}a^u \qbin{u}{v}\MM_u(x)=
a^v \MM_v(x)\, \frac{(aq^v)_{\infty}}{(ax)_{\infty}},
\]
for $v\in\N$.
The following theorem provides a $\gn$ analogue of the type
\eqref{gnseriesa}.

\begin{theorem}[$\gn$ $q$-binomial theorem I]\label{thmAL}
Let 
\begin{equation}\label{XV}
\X=x_1+\cdots+x_n \And
\V=\spec{v}_1+\cdots+\spec{v}_n,
\end{equation}
where $v\in\N^n$. Then
\begin{align}\label{qbtA}
\sum_u a^{\abs{u}} \qbin{u}{v} \MM_u(x)
&=a^{\abs{v}} \MM_v(x)
\prod_{i=1}^n \frac{(a\spec{v}_i)_{\infty}}{(ax_i)_{\infty}} \\
&=a^{\abs{v}} \MM_v(x)\:\sigma_a\biggl[\frac{\X-\V}{1-q}\biggr].
\notag
\end{align}
\end{theorem}

Note that the product in the first expression on the right
can be written in several different forms since
\begin{equation}\label{alt}
\prod_{i=1}^n (a\spec{v}_i)_{\infty}=
\prod_{i=1}^n (aq^{v_i^{+}}t^{n-i})_{\infty}=
\frac{1}{(at^{n-1})_v}
\prod_{i=1}^n (at^{n-i})_{\infty}.
\end{equation}
We also note that unlike the case $n=1$, where $v$ is readily eliminated
from the identity by shifting the summation index,
\eqref{qbtA} depends nontrivially on $v\in\N^n$ when $n>1$.

Before giving a proof we state four simple corollaries of the
$\gn$ $q$-binomial theorem.

\begin{corollary}[$\gn$ $q$-binomial theorem II]
We have
\[
\sum_u a^{\abs{u}} \MM_u(x)
=\prod_{i=1}^n \frac{(at^{n-i})_{\infty}}{(ax_i)_{\infty}}.
\]
\end{corollary}
This is of course nothing but the $v=0$ case of Theorem~\ref{thmAL}
(and corresponds to \eqref{qBT} of the introduction).
In view of \eqref{MMn1} it is, however, closer to the standard
formulation of the $q$-binomial theorem \eqref{qbt}.

A well-known finite form of the $q$-binomial theorem gives the
expansion of the $q$-shifted factorial $(a)_u$ as a power series
in $a$ \cite[Theorem 3.3]{Andrews76}
\[
(a)_u=\sum_{v=0}^u (-a)^v q^{\binom{v}{2}}\qbin{u}{v},
\qquad u\in\N.
\]
Theorem~\ref{thmAL} extends this to a $q$-shifted factorial
indexed by compositions.
\begin{corollary}
For $u\in\N^n$
\begin{equation}\label{eqav}
(a)_u=\sum_v \tau_v \, a^{\abs{v}} \qbin{u}{v}.
\end{equation}
\end{corollary}

\begin{proof}
Take Theorem~\ref{thmAL} with $(x,a,u,v)\mapsto (\spec{u},at^{1-n},v,w)$.
Using \eqref{bu} and \eqref{MMqbin} yields
\begin{equation}\label{ada}
\sum_v \tau_v \, a^{\abs{v}} \qbin{u}{v}\qbin{v}{w}
=\tau_w\,a^{\abs{w}}\frac{(a)_u}{(a)_w} \qbin{u}{w}.
\end{equation}
The identity \eqref{eqav} is the special case $w=0$.
\end{proof} 

The sum \eqref{ada} also corresponds to the $b\to 0$
limit of the $q$-Chu--Vandermonde \eqref{qCV3}. 
We stress however that the above proof is independent
of the results of Section~\ref{secE} as the proof
of Theorem~\ref{thmAL} given below only uses
vanishing properties of the Macdonald interpolation 
polynomials.

\begin{corollary}[$\sln$ Euler sum]\label{propEsum}
For $v\in \N^n$ 
\begin{equation}\label{Esum}
\sum_u \qbin{u}{v} \EE_u(x)
=\EE_v(x) \prod_{i=1}^n \frac{1}{(x_i)_{\infty}}.
\end{equation}
\end{corollary}
For $n=1$ and $v=0$ (or $u\mapsto u+v$) this is Euler's
$q$-exponential sum \cite[Equation (II.1)]{GR04}
\[
\sum_{u=0}^{\infty} \frac{x^u}{(q)_u}=\frac{1}{(x)_{\infty}}.
\]
For general $n$ it can also be found in \cite[Equation (3.30)]{BF00}.

\begin{proof}
Replace $x\mapsto x/a$ in \eqref{qbtA} and let $a$ 
tend to zero using \eqref{EEfromMM}.
\end{proof}

Our final corollary contains a generalisation of the $_1\phi_1$
summation \cite[Equation (II.5)]{GR04}.
\begin{corollary}
For $v\in\N^n$
\[
\sum_u \frac{\tau_u\,a^{\abs{u}}}{(at^{n-1})_u}
\qbin{u}{v}_{q^{-1},t^{-1}} \MM_u(x)=
\tau_v\,a^{\abs{v}} \MM_v(x)
\prod_{i=1}^n \frac{(ax_i)_{\infty}}{(at^{n-i})_{\infty}}.
\]
\end{corollary}

\begin{proof}
This results by inverting \eqref{qbtA} using \eqref{inv}.
\end{proof}

\begin{proof}[Proof of Theorem~\ref{thmAL}]
By \eqref{complete} 
the right-hand side of \eqref{qbtA} can be expanded in terms
of complete symmetric functions as
\[
\MM_v(x)\sum_{k=0}^{\infty} a^{k+\abs{v}} S_k\biggl[\frac{\X-\V}{1-q}\biggr].
\]
Comparing this with the left-hand side of \eqref{qbtA} 
and equating coefficients of $a^k$, we are left to prove that
\[
\sum_{\abs{u}=k+\abs{v}} \qbin{u}{v} \MM_u(x)=
\MM_v(x)S_k\biggl[\frac{\X-\V}{1-q}\biggr].
\]
Both sides are polynomials in $x$ of degree $k+\abs{v}$
so that it suffices to check that the above is true
for $x=\spec{w}$ where $w$ is any compositions such that 
$\abs{w}\leq k+\abs{v}$.
In other words, introducing $\W=\spec{w}_1+\cdots+\spec{w}_n$
and using \eqref{MMqbin}, we need to show that
\[
t^{(n-1)k}\sum_{\abs{u}=k+\abs{v}} \tau_u\, 
\qbin{u}{v}\qbin{w}{u}=
\tau_v \, \qbin{w}{v}
S_k\biggl[\frac{\W-\V}{1-q}\biggr]
\]
holds for all $w$ such that $\abs{w}\leq k+\abs{v}$.

Since $\abs{u}=k+\abs{v}$ it follows from \eqref{qbin0} that
the left-hand side vanishes if $\abs{w}<k+\abs{v}$. 
If on the other hand $\abs{w}=k+\abs{v}$ then, again by \eqref{qbin0}, 
only the term $u=w$ contributes to the sum. 
The previous equation can thus also be stated as
\begin{equation}\label{chi}
\qbin{w}{v} S_k\biggl[\frac{\W-\V}{1-q}\biggr]=
\frac{\tau_w}{\tau_v}\, t^{(n-1)k}\chi(\abs{w}=k+\abs{v})
\qbin{w}{v},
\end{equation}
where $\chi(\text{true})=1$ and $\chi(\text{false})=0$.

By \eqref{Knop0} both sides are identically zero if 
$v^{+}\not\subseteq w^{+}$. It is thus sufficient to show
that 
\[
S_k\biggl[\frac{\W-\V}{1-q}\biggr]=
\frac{\tau_w}{\tau_v}\, t^{(n-1)k}\chi(\abs{w}=k+\abs{v})
\]
for $v^{+}\subseteq w^{+}$ such that $\abs{w}\leq k+\abs{v}$.
Assuming these conditions we find
\[
\frac{\W-\V}{1-q}=-\Y,
\]
where $\Y=\Y_1+\cdots+\Y_n$, with 
\begin{equation}\label{Yi}
\Y_i=t^{n-i} \, \frac{q^{v_i^{+}}-q^{w_i^{+}}}{1-q}=
t^{n-i}q^{v_i^{+}}+t^{n-i}q^{v_i^{+}+1}+
\cdots+t^{n-i}q^{w_i^{+}-1}.
\end{equation}
The crucial observation is that 
\begin{subequations}
\begin{align}
\label{Ylessk}
\abs{\Y}<k&\qquad \text{if $\abs{w}<k+\abs{v}$,} \\
\abs{\Y}=k&\qquad \text{if $\abs{w}=k+\abs{v}$.}
\label{Yequalk}
\end{align}
\end{subequations}
{} From \eqref{AminB} it follows that
\[
\sigma_z[-\Y]=\prod_{y\in\Y}(1-yz)
\]
so that
\begin{subequations}\label{eq43}
\begin{align}\label{Sknul}
S_k[-\Y]&=0 && \text{if $\abs{\Y}<k$} \\
S_k[-\Y]&=(-1)^k \prod_{y\in\Y}y && \text{if $\abs{\Y}=k$.}
\label{Skk}
\end{align}
\end{subequations}
(More generally $S_k[-\Y]=(-1)^k e_k[\Y]$, with
$e_k[\Y]$ the elementary symmetric function).
Equations \eqref{Ylessk} and \eqref{Sknul} imply that
$S_k[-\Y]=0$ for $\abs{w}<k+\abs{v}$, in accordance with
\eqref{chi}. Furthermore \eqref{Yequalk} and \eqref{Skk}
imply that for $\abs{w}=k+\abs{v}$
\[
S_k[-\Y]=(-1)^k \prod_{y\in \Y} y 
=\frac{\tau_w}{\tau_v}\,t^{(n-1)k}
\]
again in agreement with \eqref{chi}.
\end{proof}

\section{A \texorpdfstring{$\gn$ $q$}{}-Gauss sum}
\label{secGauss}

To generalise the $q$-Gauss sum we 
substitute $(a,b,c,k)$ by $(q^v/x,b/a,bq^v,u)$ in \eqref{qG}, shift
the summation index and recall \eqref{MMn1}. Hence
\[
\sum_{u=0}^{\infty} a^{u-v} \frac{(b/a)_{u-v}}{(b)_u} \qbin{u}{v}
\MM_u(x) = \frac{\MM_v(x)}{(a)_v}\: 
\frac{(a,bx)_{\infty}}{(b,ax)_{\infty}}.
\]
By \eqref{homEuv} and \eqref{Euv1} this can also be written as
\[
\sum_{u=0}^{\infty}\frac{(b)_v}{(b)_u}\,\Skew{u}{v}{a}{b}\MM_u(x)
=\MM_v(x)\: \frac{(aq^v,bx)_{\infty}}{(bq^v,ax)_{\infty}},
\]
for $v\in\N$.

\begin{theorem}[$\gn$ $q$-Gauss sum I]\label{thmGauss}
With the notation of Theorem~\ref{thmAL}
\begin{align}\label{Gauss}
\sum_u \frac{(bt^{n-1})_v}{(bt^{n-1})_u}\,\Skew{u}{v}{a}{b}\MM_u(x)
&=\MM_v(x) \prod_{i=1}^n \frac{(a\spec{v}_i,bx_i)_{\infty}}
{(b\spec{v}_i,ax_i)_{\infty}}  \\
\notag 
&=\MM_v(x)\, \sigma_1\biggl[\frac{a-b}{1-q}\,(\X-\V)\biggr].
\end{align}
\end{theorem}
Note that the series on the left satisfies the $\gn$ balancing condition
\eqref{gnbalanced}. By \eqref{homEuv} and \eqref{Euvqbin} the $\gn$ $q$-Gauss 
sum simplifies to the $\gn$ $q$-binomial theorem \eqref{qbtA}
when $b$ tends to zero. 
We also remark that the $\gn$ $q$-Gauss sum may be viewed as 
a nonterminating analogue of the $q$-Chu--Vandermonde sum 
\eqref{qCV3}. Specifically, taking $(u,v,x)\mapsto (v,w,\spec{u})$
in \eqref{Gauss}, using \eqref{MMqbin} and \eqref{alt}, and finally
scaling $(a,b)\mapsto (at^{1-n},bt^{1-n})$ using \eqref{homEuv}
yields \eqref{qCV3}.

\begin{proof}[Proof of Theorem~\ref{thmGauss}]
If we take Corollary~\ref{corcon}, multiply both sides
by $b^{\abs{u}} \MM_u(x)$ and sum over $u$, we obtain
\[
\sum_v a^{-\abs{v}}\Skew{v}{w}{1}{a}\sum_u (ab)^u\qbin{u}{v}\MM_u(x)=
\sum_u b^{\abs{u}}\qbin{u}{w}\MM_u(x).
\]
Both the $u$-sums can be performed by the $\gn$ $q$-binomial theorem 
\eqref{qbtA} so that
\[
\sum_v b^{\abs{v}}\Skew{v}{w}{1}{a}\MM_v(x)
\prod_{i=1}^n \frac{(ab\spec{v}_i)_{\infty}}{(abx_i)_{\infty}}=
b^{\abs{w}}\MM_w(x)
\prod_{i=1}^n \frac{(b\spec{w}_i)_{\infty}}{(bx_i)_{\infty}}.
\]
By \eqref{bu} and \eqref{homEuv} this can be rewritten as
\[
\sum_v \frac{(abt^{n-1})_w}{(abt^{n-1})_v}\,\Skew{v}{w}{b}{ab}\MM_v(x)
=\MM_w(x) \prod_{i=1}^n \frac{(b\spec{w}_i,abx_i)_{\infty}}
{(ab\spec{w}_i,bx_i)_{\infty}},
\]
which is \eqref{Gauss} with $(a,b,u,v)\mapsto (b,ab,v,w)$.
\end{proof}

Again we give a number of simple corollaries.

\begin{corollary}[$\gn$ $q$-Gauss sum II]\label{corGauss}
We have
\begin{equation}\label{GaussS}
\sum_u a^{\abs{u}} \frac{(b)_u}{(abt^{n-1})_u}\, \MM_u(x)
=\prod_{i=1}^n \frac{(at^{n-i},abx_i)_{\infty}}{(abt^{n-i},ax_i)_{\infty}}.
\end{equation}
\end{corollary}

\begin{proof}
By \eqref{homEuv}, \eqref{bu} and \eqref{alt} this follows from
the $(v,b)\mapsto(0,ab)$ case of \eqref{Gauss}.
\end{proof}

\begin{corollary}[$\sln$ $q$-binomial theorem I]\label{corqbtI}
We have
\begin{equation}\label{qbte2}
\sum_u \Skew{u}{v}{a}{b} \EE_u(x) 
=\EE_v(x)\prod_{i=1}^n \frac{(bx_i)_{\infty}}{(ax_i)_{\infty}}
=\EE_v(x)\,\sigma_1\biggl[\frac{a-b}{1-q}\,\X \biggr].
\end{equation}
\end{corollary}

For $n=1$ and $v=0$ this is \eqref{qbt} with $(a,x)\mapsto (b/a,ax)$.

\begin{proof}
In \eqref{Gauss} we scale $(a,b,x)\mapsto (ac,bc,x/c)$.
By \eqref{homEuv} this results in
\[
\sum_u c^{\abs{u}}
\frac{(bct^{n-1})_v}{(bct^{n-1})_u}\,\Skew{u}{v}{a}{b}\MM_u(x/c)
=c^{\abs{v}}\MM_v(x/c) \prod_{i=1}^n \frac{(ac\spec{v}_i,bx_i)_{\infty}}
{(bc\spec{v}_i,ax_i)_{\infty}}.
\]
Taking the $c\to 0$ limit using \eqref{EEfromMM} yields \eqref{qbte2}.
\end{proof}

For later reference we also state the $(a,b,v)\mapsto(1,a,0)$ case
of \eqref{qbte2}.

\begin{corollary}[$\sln$ $q$-binomial theorem II]\label{corqbtII}
We have
\begin{equation}\label{qbte}
\sum_u (a)_u \EE_u(x) 
=\prod_{i=1}^n \frac{(ax_i)_{\infty}}{(x_i)_{\infty}}.
\end{equation}
\end{corollary}

This is the nonsymmetric analogue of the Kaneko--Macdonald $q$-binomial
theorem for symmetric Macdonald polynomials \cite{Kaneko96,Macdonald}
\begin{equation}\label{qbtp}
\sum_{\la} (a)_{\la} \PP_{\la}(x) 
=\prod_{i=1}^n \frac{(ax_i)_{\infty}}{(x_i)_{\infty}}.
\end{equation}
It is in fact easily shown that \eqref{qbte} and \eqref{qbtp}
are equivalent:
\[
\sum_{\la} (a)_{\la} \PP_{\la}(x)=
\sum_{\la} \sum_{u^{+}=\la} (a)_{u^{+}} \EE_u(x)=
\sum_u (a)_u \EE_u(x).
\]

\begin{corollary}
Proposition~\ref{propEg} is true.
\end{corollary}

\begin{proof}
If we multiply \eqref{Eg} by $\EE_u(x)$ and sum over $u$ using 
\eqref{struc} we obtain
\[
\sum_u\Skew{u}{v}{a}{b}\EE_u(x)=
\EE_v(x)\sum_w (b/a)_w \EE_u(ax)
=\EE_v(x/a)\prod_{i=1}^n \frac{(bx_i)_{\infty}}{(ax_i)_{\infty}}.
\]
Since this is \eqref{qbte2} the proof is complete.
\end{proof}

Let 
\[
\QQ_{\la/\mu}\biggl[\frac{a-b}{1-t}\biggr] 
:=\frac{\PP_{\la}(\spec{0})}{\PP_{\mu}(\spec{0})}\,
\PP_{\la/\mu}\biggl[\frac{a-b}{1-t}\biggr] 
\]
and $b_0b_1v_0v_1w_0w_1=qt$. Then the following skew Cauchy-type
identity is implied by \cite[Corollary 3.8]{Rains08}:
\begin{multline*}
\frac{1}{Z}
\sum_{\la} q^{\abs{\la}} \frac{(v_0,v_1)_{\la}}{(t^n,q/w_0,q/w_1)_{\la}}\,
\PP_{\la/\mu}\biggl[\frac{1-b_0}{1-t}\biggr] 
\QQ_{\la/\nu}\biggl[\frac{1-b_1}{1-t}\biggr] \\
=\frac{(q/b_0)^{\abs{\mu}} (v_0,v_1)_{\mu}}
     {(t^n,q/b_0w_0,q/b_0w_1)_{\mu}}\,
\frac{(q/b_1)^{\abs{\nu}} (v_0,v_1)_{\nu}}
     {(t^n,q/b_1w_0,q/b_1w_1)_{\nu}} \qquad\qquad\qquad\qquad\qquad\\
\times \sum_{\la} (b_0b_1/q)^{\abs{\la}}
\frac{(t^n,q/b_0b_1w_0,q/b_0b_1w_1)_{\la}}{(v_0,v_1)_{\la}}\,
\PP_{\nu/\la}\biggl[\frac{1-b_0}{1-t}\biggr]
\QQ_{\mu/\la}\biggl[\frac{1-b_1}{1-t}\biggr],
\end{multline*}
provided the sum on the sum on the left terminates.
The normalisation $Z$ is given by the sum on the left for $\mu=\nu=0$.
Choosing
\[
(b_0,b_1,v_0,v_1,w_0,w_1)\mapsto (b/a,d/c,t^n,q^{-k},qt^{1-n}/bd,q^k ac)
\]
with $k$ a nonnegative integer and taking the large $k$ limit yields
\begin{multline}\label{Zp}
\frac{1}{Z'}
\sum_{\la}\frac{1}{(bdt^{n-1})_{\la}}\,
\PP_{\la/\mu}\biggl[\frac{a-b}{1-t}\biggr] 
\QQ_{\la/\nu}\biggl[\frac{c-d}{1-t}\biggr] \\
=\frac{1}{(adt^{n-1})_{\mu}(bct^{n-1})_{\nu}} 
\sum_{\la} (act^{n-1})_{\la}\,
\PP_{\nu/\la}\biggl[\frac{a-b}{1-t}\biggr]
\QQ_{\mu/\la}\biggl[\frac{c-d}{1-t}\biggr],
\end{multline}
where $Z'$ again denotes the sum on the left for $\mu=\nu=0$.
Explicitly,
\[
Z'=\sum_{\la}\frac{(b/a,d/c)_{\la}}{(bdt^{n-1})_{\la}}\,
\PP_{\la}(ac\spec{0}) \\
=\prod_{i=1}^n \frac{(adt^{n-i},bct^{n-i})_{\infty}}
{(act^{n-i},bdt^{n-i})_{\infty}},
\]
where the second equality follows from Kaneko's $q$-Gauss
sum for Macdonald polynomials \cite[Proposition 5.4]{Kaneko96}.
Consequently, \eqref{Zp} may also be written as
\begin{multline*}
\sum_{\la}\biggl(\:\prod_{i=1}^n (bd\spec{\la}_i)_{\infty}\biggr)
\PP_{\la/\mu}\biggl[\frac{a-b}{1-t}\biggr] 
\QQ_{\la/\nu}\biggl[\frac{c-d}{1-t}\biggr] \\
=\sum_{\la}
\biggl(\:\prod_{i=1}^n \frac{(ad\spec{\mu}_i,bc\spec{\nu}_i)_{\infty}}
{(ac\spec{\la}_i)_{\infty}}\biggr)
\PP_{\nu/\la}\biggl[\frac{a-b}{1-t}\biggr]
\QQ_{\mu/\la}\biggl[\frac{c-d}{1-t}\biggr].
\end{multline*}
As our final corollary of the $\gn$ $q$-Gauss sum we will show that this
has a nonsymmetric analogue.
Let 
\[
\SkewF{u}{v}{a}{b}:=\frac{\EE_u(\spec{0})}{\EE_v(\spec{0})}\,
\Skew{u}{v}{a}{b}
\]
\begin{corollary}
For $v,w\in\N^n$ and $\abs{ac},\abs{ad},\abs{bc},\abs{bd}$ sufficiently
small so that the sum on the left converges
\begin{multline*}
\sum_u \biggl(\:\prod_{i=1}^n (bd\spec{u}_i)_{\infty}\biggr)
\Skew{u}{v}{a}{b} \SkewF{u}{w}{c}{d} \\
=\sum_u \biggl(\:\prod_{i=1}^n 
\frac{(ad\spec{v}_i,bc\spec{w}_i)_{\infty}}{(ac\spec{u}_i)_{\infty}}\biggr)
\Skew{w}{u}{a}{b} \SkewF{v}{u}{c}{d}.
\end{multline*}
\end{corollary}

\begin{proof}
Application of \eqref{cdef} to both sides of the $\gn$ $q$-Gauss sum
results in
\[
\sum_{u,w} \frac{(bt^{n-1})_v}{(bt^{n-1})_u}\,\Skew{u}{v}{a}{b}
c_{uw}(1,c) \MM_w(cx)
=\sum_u c_{vu}(1,c) \MM_u(cx)
\prod_{i=1}^n \frac{(a\spec{v}_i,bx_i)_{\infty}}
{(b\spec{v}_i,ax_i)_{\infty}}.
\]
The right-hand side can again be transformed by the $q$-Gauss sum 
\eqref{Gauss} with $(a,b,x,v)\mapsto (a/c,b/c,cx,w)$ so that
\begin{multline*}
\sum_{u,w}
\biggl(\:\prod_{i=1}^n (b\spec{u}_i)_{\infty}\biggr)
\Skew{u}{v}{a}{b} c_{uw}(1,c) \MM_w(cx) \\
=\sum_{u,w} 
\biggl(\:\prod_{i=1}^n \frac{(a\spec{v}_i,b\spec{w}_i/c)_{\infty}}
{(a\spec{u}_i/c)_{\infty}}\biggr)
c_{vu}(1,c) 
\Skew{w}{u}{a/c}{b/c} \MM_w(cx),
\end{multline*}
where we have also used \eqref{bu}. 
Next we equate coefficients of $\MM_w(cx)$ and
substitute $(a,b,c)\mapsto (ad,bd,d/c)$.
After carrying out some simplifications
using \eqref{Skewdef} and \eqref{homEuv} the
result follows.
\end{proof}

\section{A \texorpdfstring{$\gn$ $q$}{}-Kummer--Thomae--Whipple formula}
\label{secKTW}

By the substitutions 
\[
(a,b,d,e,f)\mapsto (1/x,bq^v,dq^v,eq^v,a)
\]
and the use of \eqref{homEuv} and \eqref{Euv1}, the
Kummer--Thomae--Whipple formula \eqref{KTW1} can be written in the form
\[
\sum_{u=0}^{\infty}\frac{(b)_u}{(d,e)_u}\,\Skew{u}{v}{a}{ac}\MM_u(x)
=\frac{(b)_v}{(d/c)_v}\,
\frac{(a,ex)_{\infty}}{(e,ax)_{\infty}}
\sum_{u=0}^{\infty} \frac{(d/c)_u}{(d,a)_u}\, 
\Skew{u}{v}{e}{ac}\MM_u(x).
\]

\begin{theorem}[$\gn$ $q$-Kummer--Thomae--Whipple formula]\label{thmKTW}
For $v\in\N^n$ and $abc=de$
\begin{multline*}
\sum_u \frac{(b)_u}{(d,et^{n-1})_u}\,\Skew{u}{v}{a}{ac}\MM_u(x) \\
=\frac{(b)_v}{(d/c)_v}\,
\biggl(\:\prod_{i=1}^n \frac{(at^{n-i},ex_i)_{\infty}}
{(et^{n-i},ax_i)_{\infty}} \biggr) 
\sum_u \frac{(d/c)_u}{(d,at^{n-1})_u}\,\Skew{u}{v}{e}{ac}\MM_u(x) .
\end{multline*}
\end{theorem}
Note that the condition $abc=de$ corresponds to the $\gn$ balancing condition
\eqref{gnbalanced} and that the transformation may alternatively be written
as
\begin{multline*}
\sum_u \frac{(b)_u(d,et^{n-1})_v}{(b)_v(d,et^{n-1})_u}\, 
\Skew{u}{v}{a}{ac}\MM_u(x) \\
=\biggl(\:\prod_{i=1}^n \frac{(a\spec{v}_i,ex_i)_{\infty}}
{(e\spec{v}_i,ax_i)_{\infty}} \biggr) 
\sum_u \frac{(d/c)_u(d,at^{n-1})_v}{(d/c)_v(d,at^{n-1})_u}\, 
\Skew{u}{v}{e}{ac}\MM_u(x),
\end{multline*}
with the product on the right corresponding to
$\sigma_1\bigl[(a-e)(\X-\V)/(1-q)\bigr]$ in the notation 
of Theorem~\ref{thmAL}.

The above theorem generalises a number of earlier results.
For example, recalling \eqref{deltaE} and taking the limit $e\to ac$ 
(so that $d\to b$), and finally replacing $c\mapsto b/a$ we obtain the 
$\gn$ $q$-Gauss sum of Theorem~\ref{thmGauss}.
Furthermore, making the substitutions 
$(a,e,u,v,x)\mapsto(at^{1-n},et^{1-n},v,w,\spec{u})$,
and using \eqref{bu}, \eqref{MMqbin} and \eqref{homEuv},
we arrive at \eqref{Searsn2}.

\begin{proof}[Proof of Theorem~\ref{thmKTW}]
If we replace $(a,c,d)\mapsto(d/c,d,ac/d)$ in \eqref{qps}
and define $e=abc/d$ we obtain
\begin{equation}\label{qps2}
\sum_v \frac{(d/c)_v}{(d)_v}\,\Skew{u}{v}{a}{e}\Skew{v}{w}{e}{ac}
=\frac{(d/c)_w(b)_u}{(b)_w(d)_u}\,\Skew{u}{w}{a}{ac}.
\end{equation}
Now denote the left-hand side of Theorem~\ref{thmKTW} with
$v$ replaced by $w$ as $\text{LHS}$. By \eqref{qps2},
\begin{align*}
\text{LHS}&=
\sum_u \frac{(b)_u}{(d,et^{n-1})_u}\,\Skew{u}{w}{a}{ac}\MM_u(x) \\
&=\frac{(b)_w}{(d/c)_w} \sum_{u,v} \frac{(d/c)_v}{(d)_v}\,
\frac{1}{(et^{n-1})_u}\,\Skew{u}{v}{a}{e}\Skew{v}{w}{e}{ac}\MM_u(x).
\end{align*}
The sum over $u$ can now be performed by the $\gn$ $q$-Gauss 
sum \eqref{Gauss} with $b\mapsto e$. Hence
\[
\text{LHS}=\frac{(b)_w}{(d/c)_w}
\biggl(\:\prod_{i=1}^n \frac{(ex_i)_{\infty}}{(ax_i)_{\infty}} \biggr) 
\sum_v \biggl(\:\prod_{i=1}^n \frac{(a\spec{v}_i)_{\infty}}
{(e\spec{v}_i)_{\infty}}\biggr)
 \frac{(d/c)_v}{(d)_v}\,\Skew{v}{w}{e}{ac}\MM_v(x).
\]
Finally using \eqref{alt} this yields the right-hand side of
the theorem with $v$ replaced by $w$.
\end{proof}

A number of new results follow from the $\gn$ 
$q$-Kummer--Thomae--Whipple formula.

\begin{corollary}[$\gn$ Heine transformation]\label{thmOW}
For $v\in\N^n$
\begin{multline*}
\sum_u a^{\abs{u}}\frac{(b)_u}{(ct^{n-1})_u}\qbin{u}{v}\MM_u(x) \\
=\Bigl(\frac{a}{c}\Bigr)^{\abs{v}} \frac{(b)_v}{(ab/c)_v}\,
\biggl(\:\prod_{i=1}^n \frac{(at^{n-i},cx_i)_{\infty}}
{(ct^{n-i},ax_i)_{\infty}} \biggr)
\sum_u c^{\abs{u}} \frac{(ab/c)_u}{(at^{n-1})_u}
\qbin{u}{v} \MM_u(x).
\end{multline*}
\end{corollary}
An equivalent way to state the above transformation is 
\begin{multline*}
\sum_u a^{\abs{u}-\abs{v}} \frac{(b)_u(ct^{n-1})_v}{(b)_v(ct^{n-1})_u}
\qbin{u}{v}\MM_u(x) \\
=\biggl(\:\prod_{i=1}^n \frac{(a\spec{v}_i,cx_i)_{\infty}}
{(c\spec{v}_i,ax_i)_{\infty}} \biggr)
\sum_u c^{\abs{u}-\abs{v}} \frac{(ab/c)_u(at^{n-1})_v}{(ab/c)_v(at^{n-1})_u}
\qbin{u}{v}\MM_u(x).
\end{multline*}
For $n=1$ the $\gn$ Heine transformation simplifies to \eqref{Heine1}
with $(a,b,c,x)$ replaced by $(q^v/x,bq^v,cq^v,ax)$.

\begin{proof}
Let $(c,d)\to(0,0)$ in Theorem~\ref{thmKTW} such that $d/c=ab/e$.
By \eqref{homEuv} and \eqref{Euvqbin} the $\gn$ Heine transformation 
with $c\mapsto e$ follows.
\end{proof}

\begin{corollary}[$\sln$ $q$-Euler transformation]\label{slnEulertrafo}
For $v\in\N^n$
\[
\sum_u \frac{(a)_u}{(c)_u}\,\Skew{u}{v}{b}{c}\EE_u(x)
=\frac{(a)_v}{(b)_v}\, 
\biggl(\:\prod_{i=1}^n \frac{(ax_i)_{\infty}}{(bx_i)_{\infty}} \biggr) 
\sum_u \frac{(b)_u}{(c)_u}\,\Skew{u}{v}{a}{c} \EE_u(x).
\]
\end{corollary}
For $n=1$ this is \eqref{Heine2} with $(a,b,c,x)\mapsto(aq^v,c/b,cq^v,bx)$,
and for $v=0$ it is the nonsymmetric analogue of \cite[Proposition 3.1]{BF99}.
The $\sln$ Euler transformation is easily seen to be equivalent to the
$q$-Pfaff--Saalsch\"utz sum \eqref{qPS}. Indeed, if we take the latter,
multiply both sides by $\EE_u(x)$ and then sum over $u$ using \eqref{qbte2}
we obtain the former. 
We also note that for $c=0$ the $\sln$ $q$-Euler transformation 
becomes
\[
\sum_u b^{\abs{u}} (a)_u \qbin{u}{v}\EE_u(x)
=\Bigl(\frac{b}{a}\Bigr)^{\abs{v}} \frac{(a)_v}{(b)_v}\,
\biggl(\:\prod_{i=1}^n \frac{(ax_i)_{\infty}}{(bx_i)_{\infty}} \biggr) 
\sum_u a^{\abs{u}} (b)_u \qbin{u}{v} \EE_u(x).
\]
When $n=1$ this is trivial since both sides are summable by the $q$-binomial
theorem \eqref{qbt}. Curiously, for $n>1$ it no longer appears possible to
explicitly perform the sums in closed form.

\begin{proof}
Corollary~\ref{slnEulertrafo} follows from Theorem~\ref{thmKTW}
by replacing $(x,e)\mapsto (x/a,ae)$, then taking the  
limit $a\to 0$ using \eqref{EEfromMM}, and finally
making the substitutions
\[
(x,b,c,d,e)\mapsto (bx,a,c/b,c,a/b). \qedhere
\]
\end{proof}

\subsection*{Acknowledgements}
Part of this work was carried out at MSRI during 
the special programme \emph{Combinatorial Representation Theory}.
AL and SOW thank MSRI for hospitality and financial support.

\bibliographystyle{amsplain}

\end{document}